\documentclass[smallcondensed,nospthms]{svjour3}
\usepackage[paperwidth=15.51cm,
paperheight=23.51cm,left=1.75cm,
right=1.75cm,top=2cm,bottom=2cm]{geometry}

\usepackage{fix-cm}
\usepackage{amsmath,amsthm,amssymb,amsfonts}
\usepackage{marvosym}
\usepackage{bm}
\usepackage{color}
\usepackage{rotating}
\usepackage{subcaption}
\usepackage{psfrag,graphicx,pstool,csquotes}

\usepackage{algorithm}
\usepackage{algpseudocode}
\usepackage{lineno}

\theoremstyle{plain}
\newtheorem{thm}{Theorem}[section]
\newtheorem{lem}[thm]{Lemma}
\newtheorem{prop}[thm]{Proposition}

\theoremstyle{definition}
\newtheorem{defn}{Definition}[section]
\newtheorem{exmp}{Example}[section] 

\theoremstyle{remark}
\newtheorem{rem}[thm]{Remark}

\everymath{\displaystyle}
\makeatletter
\renewcommand*\env@matrix[1][*\c@MaxMatrixCols c]{%
	\hskip -\arraycolsep
	\let\@ifnextchar\new@ifnextchar
	\array{#1}}
\makeatother

\makeatletter
\newcommand*\bigcdot{\mathpalette\bigcdot@{.5}}
\newcommand*\bigcdot@[2]{\mathbin{\vcenter{\hbox{\scalebox{#2}{$\m@th#1\bullet$}}}}}
\makeatother

\journalname{}

\makeatletter
\providecommand{\doi}[1]{%
  \begingroup
    \let\bibinfo\@secondoftwo
    \urlstyle{rm}%
    \href{http://dx.doi.org/#1}{%
      doi:\discretionary{}{}{}%
      \nolinkurl{#1}%
    }%
  \endgroup
}
\makeatother

\begin{document}
\newcommand{\titletext}{ 
From big q-Jacobi and Chebyshev polynomials to exponential-reproducing subdivision: new identities
}
\newcommand{\titlerunningtext}{
Big q-Jacobi, Chebyshev polynomials and exponential-reproducing subdivision
}

\newcommand{\abstracttext}{
In this paper we derive
new identities satisfied by Chebyshev polynomials of the first kind and big q-Jacobi polynomials.
An immediate benefit of the derived identities is the achievement of closed-form expressions for the Laurent  polynomials that identify minimum-support interpolating 
subdivision schemes reproducing finite sets of integer
powers of exponentials.
}

\newcommand{\separ}{\and}

\newcommand{\keytext}{
Big q-Jacobi polynomial \separ Chebyshev polynomial \separ Exponential polynomial \separ
Non-stationary interpolatory subdivision    
}

\newcommand{\MSCtext}{
33C45 \separ 42C05 \separ 65D05 \separ 65D17 
}

\title{\titletext}
\titlerunning{\titlerunningtext} 

\author{
Leonard Peter Bos \and
Lucia Romani \and
Alberto Viscardi
}
\authorrunning{L.P.~Bos, L.~Romani, A.~Viscardi}  

\institute{
L.P.~Bos \at
Dipartimento di Informatica, Universit\`{a} di Verona,
Str. le Grazie 15, 37134 Verona, Italy\\
\email{leonardpeter.bos@univr.it}
\and
L.~Romani \at
Dipartimento di Matematica, Alma Mater Studiorum Universit\`{a} di Bologna,
Piazza di Porta San Donato 5, 40126 Bologna, Italy\\
\email{lucia.romani@unibo.it}
\and
A.~Viscardi \at
Dipartimento di Matematica ``Tullio Levi-Civita'', Universit\`{a} di Padova,
Via Trieste 63, 35131 Padova, Italy\\
\email{alberto.viscardi@unipd.it} 
}

\date{ }

\maketitle

\begin{abstract}
\abstracttext
\keywords{\keytext}
\subclass{\MSCtext}
\end{abstract}

{ 
\section{Main results and overview of the paper}
Big q-Jacobi polynomials, introduced by Andrews and Askey \cite{AA85}, are a family of q-orthogonal polynomials that generalize classical Jacobi polynomials within the theory of q-hypergeometric functions.
Precisely, for $n\in \mathbb{N}_0=\mathbb{N} \cup \{0\}$ and real parameters $\alpha,\beta>-1$, the classical Jacobi polynomial $P_n^{(\alpha,\beta)}(x)$ is a degree-$n$ polynomial defined by the Gauss hypergeometric function $_2F_1$ (see, e.g., \cite[eq. (9.8.1)]{KSbook}) as
$$
P_n^{(\alpha,\beta)}(x)= \frac{(\alpha + 1)_n}{n!} \ {}_2F_1 \left( 
-n, \ n+\alpha+\beta+1; \alpha+1; \frac{1-x}{2}
\right),\quad -1<x<1,
$$
where $(\alpha+1)_n$ is a Pochhammer symbol,
while the big q-Jacobi polynomial $P_n(x;a,b,c;q)$ is defined by the q-hypergeometric function $_3\phi_2$ (see, e.g., \cite[eq. (14.5.1)]{KSbook}) as 
$$
P_n(x;a,b,c;q)={}_3\phi_2 \left[ 
\begin{array}{c}
q^{-n}, \ a b q^{n+1}, \ x\\
a q, \ c q   
\end{array}
; q, q
\right],\qquad 0<q<1. 
$$
The well-known identity (see, e.g., \cite[eq. (14.5.17)]{KSbook})
$$
\lim_{q \rightarrow 1^{-}} P_n(x;q^{\alpha}, q^{\beta}, -q^{\gamma}; q)\;=\;\frac{n!}{(\alpha + 1)_n} P_n^{(\alpha,\beta)}(x),
$$
satisfied for any arbitrary real $\gamma$, links classical Jacobi polynomials to big q-Jacobi polynomials and, 
when $\alpha=\beta=-1/2$, 
allows us to write that, 
for any arbitrary real $\gamma$,
\begin{equation}\label{eq:Jacobi_as_bigqlimit}
\lim_{q \rightarrow 1^{-}} P_n(x;q^{-\frac{1}{2}}, q^{-\frac{1}{2}}, -q^{\gamma}; q)\;=\; \frac{n!}{\left(\frac{1}{2}\right)_n}P_n^{(-\frac{1}{2},-\frac{1}{2})}(x).
\end{equation}
Thus, recalling that (see, e.g., \cite[eq. (9.8.35)]{KSbook}) the degree-$n$ Chebyshev polynomial $T_n(x)$ satisfies
\begin{equation}\label{eq:Chebyshev_linked_Jacobi}
T_n(x)\;=\;{}_2F_1 \Big(-n, n; \frac{1}{2}; \frac{1-x}{2}  \Big)\;=\;\frac{n!}{\left(\frac{1}{2}\right)_n} P_n^{(-\frac{1}{2},-\frac{1}{2})}(x),
\end{equation} 
equation \eqref{eq:Jacobi_as_bigqlimit} implies also that
$$
\lim_{q \rightarrow 1^{-}} P_n(x;q^{-\frac{1}{2}}, q^{-\frac{1}{2}}, -q^{\gamma}; q)=T_n(x).
$$
However, to the best of our knowledge, no one has ever proved before that, for 
$x=(t+t^{-1})/2$, the identity
\begin{equation}\label{eq:main_result1}
\left( T_n \left(\frac{t+t^{-1}}{2} \right) \right)^{-1}=P_n\left(t;-t^{-1}, -t^{-1}, -1; t^2\right)
\end{equation}  
holds true,
which means that, for any arbitrary real $\gamma$,
$$
\lim_{q \rightarrow 1^{-}} P_n\left(\frac{t+t^{-1}}{2};q^{-\frac{1}{2}}, q^{-\frac{1}{2}}, -q^{\gamma}; q\right)\;=\;\left(\;P_n\left(t;-t^{-1}, -t^{-1}, -1; t^2\right)\;\right)^{-1}.
$$
By combining \eqref{eq:main_result1} and \eqref{eq:Chebyshev_linked_Jacobi} with the above specified expression of $x$, the new identity 
$$
P_n^{(-\frac{1}{2},-\frac{1}{2})}\left( \frac{t + t^{-1}}{2} \right) \;=\; \frac{\left(\frac{1}{2}\right)_n}{n! \; P_n\left(t;-t^{-1}, -t^{-1}, -1; t^2\right)}
$$
is also obtained.\\
Another remarkable consequence of \eqref{eq:main_result1} is that, for all $n\in \mathbb{N}$ and for 
$$
C_{\ell,i}(x)\;=\;\frac{ \big( T_{\ell}(x) - T_{\ell+1}(x) \big) \ \big( T_{\ell}(x)+1 \big)}{\big( T_{i}(x) - T_{i+1}(x) \big) \, - \, \big( T_{\ell}(x) - T_{\ell+1}(x) \big)}, \qquad 
\begin{array}{l}
    i=1,\ldots,n,\\
    \ell=0,\ldots,i-1,
\end{array}
$$
the more sophisticated identity
\begin{equation} \label{eq:main_result2}
\sum_{i=1}^n \frac{2^i}{T_i(x)+1} \; \prod_{\ell=0}^{i-1}  C_{\ell,i}(x) \, \frac{\Big(T_{\ell}(x)\Big)^2 - \Big(T_{n}(x)\Big)^2}{\Big( T_{\ell}(x) + 1 \Big)^2}  \;=\; \frac{1-T_n(x)}{2T_n(x)}  
\end{equation}
holds true as well.
This last result provides the keystone to arrive at the closed-form expression for the symbol $m_{2n+2,k}(z)$, $k\in\mathbb{N}_0$, of the $k$-level Laurent polynomial of the minimum-support interpolating subdivision scheme 
generating the function space 
\begin{equation} \label{eq:V_space}
V_{2n+2,\theta}\;=\;{\rm span}\left\{\;1,\;(\bigcdot), \;  \{e^{\pm \mathrm{i} j\theta \, \bigcdot}\}_{j=1}^{n}\;\right\}\;\cap\;\mathcal{C}^\infty(\mathbb{R},\mathbb{R}),
\end{equation}
with $\theta \in [0, {\sigma_{\star}(n)}[\; \cup\; \mathrm{i} \mathbb{R}^+$ and $\sigma_{\star}(n)\;\in\;]0,2\pi[$ the so-called critical length of the space (see, e.g., \cite{BCM19,BCM20,CMP04}), i.e.,  
\begin{equation} \label{eq:main_result3}
m_{2n+2,k}(z)\;=\;
2 a_{0,k}(z) + 
2 \big(2 a_{0,k}(z)-1\big) \, 
\sum_{i=1}^n c_i(v_k) \,  
\prod_{\ell=0}^{i-1} a_{\ell,k}(-z) \,  a_{\ell,k}(z),  
\end{equation}
with
$$ 
c_i(v_k)=\frac{2^i \displaystyle \prod_{\ell=0}^{i-1} C_{\ell,i}(v_{k})}{T_i(v_{k}) +1}, \quad
    a_{\ell,k}(z)=\frac{z + 2 T_{\ell}\big( v_{k} \big) +z^{-1}}{2 \Big( T_{\ell}\big( v_{k} \big)+1 \Big)},  \quad  v_k=\cos \left( \frac{\theta}{2^{k+1}} \right).   
$$ 
The function space $V_{2n+2,\theta}$, mixing linear polynomials with either hyperbolic or trigonometric functions (depending on the choice of $\theta$), and some of its subspaces have been extensively considered in literature, e.g., in \cite{BCR07,CR2010,MP10,Romani09}.
The interpolatory subdivision scheme with symbols $\{m_{2n+2,k}(z)\}_{k \in \mathbb{N}_0}$ turns out to be the non-stationary version of the $(2n+2)$-point Dubuc-Deslauriers interpolating scheme reproducing algebraic polynomials \cite{DD}, for which a closed-form expression was first derived in  \cite{DS2006} by using a Fourier  
analysis approach and later in \cite{DHSS} by means of an algebraic approach. 
Although in \cite{DLL03} it was already shown that, as $k \rightarrow +\infty$, the 
$k$-level Laurent polynomial of any non-stationary interpolatory subdivision scheme reproducing a $(2n+2)$-dimensional space of exponential polynomials converges to the symbol of the $(2n+2)$- point Dubuc–Deslauriers interpolatory scheme and the two share the same smoothness, no closed-form expression for the symbols of such a minimum-support family of non-stationary interpolatory schemes was there presented.\\

The remainder of the paper is organized as follows. In Section \ref{sec:1a} we start introducing the notation and recalling the preliminary definitions that in Section \ref{sec:1b} allow us to present our new results on Chebyshev polynomials of the first kind. In particular, equations \eqref{eq:main_result1} and \eqref{eq:main_result2} are proven in Proposition \ref{prop:1_su_Tn} and Proposition \ref{prop:identity_Tn}, respectively.
Then, after describing in Section \ref{sec:0} the state-of-the-art on minimum-support interpolating subdivision schemes reproducing finite sets of integer powers of exponentials, in Section \ref{sec:3} we exploit the identities derived in Section \ref{sec:1b} to construct the closed-form expressions of their symbols as in \eqref{eq:main_result3} (see Theorem \ref{theo_main}). Section \ref{sec:closure} concludes the
paper with some final remarks.
}

\section{Preliminaries and notation}\label{sec:1a}

\begin{defn}
Let $n \geq 0$ be an integer.
The Pochhammer symbol $(p)_n$ is nothing but a shifted factorial defined by
$$
(p)_n=\left\{
\begin{array}{ll}
1 & \hbox{if} \quad n=0,\\
\displaystyle \prod_{j=0}^{n-1} (p+j) & \hbox{if} \quad n>0.
\end{array}
\right.
$$
\end{defn}
Special values include $(1)_n=n!$ and $\left( \frac{1}{2} \right)_n=\frac{(2n-1)!!}{2^n}=\frac{(2n)!}{(2^n)^2 n!}$.\\ 

The Pochhammer symbol plays a key role in the definition of hypergeometric functions. 

\begin{defn}
The Gauss hypergeometric function ${}_2F_1(a,b;c;z)$ is defined for $|z|<1$ by the power series  
$$
{}_2F_1(a,b;c;z)=\sum_{n=0}^{\infty} \frac{(a)_n \, (b)_n}{(c)_n} \ \frac{z^n}{n!}. 
$$
The series is undefined (or infinite) if $c$ is a non-positive integer.
The series terminates if either $a$ or $b$ is a non-positive integer, in which case the function reduces to a polynomial of the form 
$$
{}_2F_1(-m,b;c;z)=\sum_{n=0}^{m} (-1)^n \, \binom{m}{n} \,  \frac{(b)_n}{(c)_n} \, z^n. 
$$
\end{defn}

Many of the common mathematical functions can be expressed in terms of the Gauss hypergeometric function and several orthogonal polynomials, including Jacobi polynomials and their special cases, can be written in terms of hypergeometric functions.
In particular (cf. \cite[page 2]{GasperRahman}, Equation (1.2.3)), the degree-$n$ Chebyshev polynomial  satisfies
\begin{equation}\label{eq:Tn_vs_2F1}
{}_2F_1 \Big(-n, n; \frac{1}{2}; \frac{1-x}{2} \Big)=T_n(x)=\frac{(2^n n!)^2}{(2n)!} \, P_n^{(-\frac{1}{2},-\frac{1}{2})}(x),
\end{equation}
where the last equality comes from the well-known relationship between Chebyshev and Jacobi polynomials. 

\smallskip 
The Gauss hypergeometric function proved to be a powerful tool also in the context of subdivision schemes \cite{RVJima} to provide the limit-point stencils for a class of non-stationary subdivision schemes generalizing the ones in \cite{GP00}.

\begin{defn}
The q-analogue of the Pochhammer symbol is the q-Pochhammer symbol $(a;q)_n$ defined by
$$
(a;q)_n=\left\{
\begin{array}{ll}
1 & \hbox{if} \quad n=0,\\
\displaystyle \prod_{j=0}^{n-1} (1-aq^j) & \hbox{if} \quad n>0.
\end{array}
\right.
$$
\end{defn}

The q-Pochhammer symbol is linked to the Pochhammer symbol by the fact that
$$
\lim_{q \rightarrow 1^{-}} \frac{(q^p;q)_n}{(1-q)^n}=(p)_n,
$$
and is a major building block in the theory of generalized hypergeometric series: indeed it plays the role that the ordinary Pochhammer symbol plays in the theory of basic  hypergeometric series.

\begin{defn}
The q-hypergeometric function ${}_r\phi_s$ is  
$$
{}_r\phi_s \left[ 
\begin{array}{ccc}
\alpha_1, & \ldots, & \alpha_r\\
\beta_1, & \ldots, & \beta_s
\end{array}
; q, z
\right]
= \sum_{n=0}^{\infty} \frac{(\alpha_1;q)_n \, \ldots \, (\alpha_r;q)_n}{(\beta_1;q)_n \, \ldots \, (\beta_s;q)_n} \ \frac{z^n}{(q;q)_n} \ 
\left(-q^\frac{n-1}{2}\right)^{n(1+s-r)}
$$
where $q\neq 0$ when $r>s+1$.
\end{defn}

\medskip 
{  
We observe that ${}_r\phi_s$ is defined on  $]-1, 1[$ if $r = s + 1$ and on $\mathbb{R}$ if $r < s + 1$. If one of the $\alpha_i$ is equal to $q^{-n}$ for some integer $n \geq 0$, then ${}_r\phi_s$ is a degree-$n$ polynomial. 
In this case one of the $\beta_j$ can be equal to $q^{-m}$ for some integer $m \geq 0$ if $m \geq n$.\\
}
Since
$$
\lim_{q \rightarrow 1^{-}} {}_r\phi_s \left[ 
\begin{array}{ccc}
q^{\alpha_1}, & \ldots, & q^{\alpha_r}\\
q^{\beta_1}, & \ldots, & q^{\beta_s}
\end{array}
; q, (q-1)^{1+s-r} z
\right] 
\;=\;
{}_rF_s \left[ 
\begin{array}{ccc}
\alpha_1, & \ldots, & \alpha_r\\
\beta_1, & \ldots, & \beta_s
\end{array}
; z
\right]
$$
where ${}_rF_s$ is a generalized hypergeometric function defined as
$$
{}_rF_s \left[ 
\begin{array}{ccc}
\alpha_1, & \ldots, & \alpha_r\\
\beta_1, & \ldots, & \beta_s
\end{array}
; z
\right]
\;=\;
\sum_{n=0}^{\infty} \frac{(\alpha_1)_n \ldots (\alpha_r)_n}{(\beta_1)_n \ldots (\beta_s)_n} \ \frac{z^n}{n!},
$$ 
the function ${}_r\phi_s$ is considered a q-analogue of ${}_rF_s$.

One of the most important formulas for q-hypergeometric functions is the 
q–Saalsch\"{u}tz summation formula 
(cf. \cite[page 355]{GasperRahman}, Equation (II.12)), also known as Jackson's summation formula. It applies to the case $r=3$, $s=2$ with 
$$
\alpha_1= q^{-n}, 
\qquad 
\alpha_2=a, 
\qquad
\alpha_3=b, 
\qquad 
\beta_1=c,
\qquad
\beta_2=\frac{ab}{cq^{n-1}}, 
\qquad 
z=q,
$$
and reads as
\begin{equation}\label{eq:q-pfaff}
{}_3\phi_2 \left[ 
\begin{array}{c}
q^{-n}, \ a, \ b\\
c, \ \frac{ab}{cq^{n-1}}  
\end{array}
; q, q
\right] 
\;=\;
\frac{(\frac{c}{a};q)_n \ (\frac{c}{b};q)_n}{(c;q)_n \ (\frac{c}{ab}; q)_n}, \qquad n \in \mathbb{N}.
\end{equation}

The q-hypergeometric function ${}_3\phi_2$ 
plays a key role in the domain of special functions since it is exploited to define 
q-analogues of classical Jacobi polynomials (see, e.g., \cite[eq. (14.5.1)]{KSbook} or \cite{GasperRahman} for details), i.e., 
big q-Jacobi polynomials $P_n(z;a,b,c;q)$ as
\begin{equation}\label{eq:bigq-Jac}
P_n(z;a,b,c;q)
\;=\;
{}_3\phi_2 \left[ 
\begin{array}{c}
q^{-n}, \ a b q^{n+1}, \ z\\
a q, \ c q   
\end{array}
; q, q
\right].
\end{equation}
In the following section we prove a new result linking big q-Jacobi polynomials to Chebyshev polynomials.  

\section{Chebyshev polynomials of the first kind: new results}\label{sec:1b}

By introducing the change of variable
\begin{equation}\label{eq:xvst}
x\;=\;\frac{t+t^{-1}}{2},
\end{equation}
we rewrite the $n$-th Chebyshev polynomial as
$T_n(x)=(t^n + t^{-n})/2$, and prove a new result linking $\big(T_n(x) \big)^{-1}$ to the q-hypergeometric function ${}_3\phi_2$ defined in terms of $t$. We point out that the resulting identity, as well as the following results of this section, should be intended as algebraic identities since both sides could be not well-defined for some values of $t$ and $x$. 

\begin{prop}\label{prop:1_su_Tn}
$$
{}_3\phi_2 \left[ 
\begin{array}{c}
t^{-2n}, \ t^{2n}, \ t\\
-t, \ -t^2  
\end{array}
; t^2, t^2
\right] 
\;=\; 
\frac{2t^n}{1+t^{2n}}  
\;=\;
\frac{1}{T_n(x)}, \qquad n \in \mathbb{N}_0.
$$
\end{prop}
\begin{proof}
By applying formula \eqref{eq:q-pfaff} with $q=t^2$, $a=t^{2n}$, $b=t$, $c=-t$
we are able to write
$$
{}_3\phi_2 \left[ 
\begin{array}{c}
t^{-2n}, \ t^{2n}, \ t\\
-t, \ -t^2  
\end{array}
; t^2, t^2
\right] 
\;=\; \frac{(-1;t^2)_n \ (-t^{-(2n-1)};t^2)_n}{(-t;t^2)_n \ (-t^{-2n};t^2)_n}.
$$
Now, taking into account that
$$
\begin{array}{rcl}
    (-1;t^2)_n &=& \displaystyle \prod_{j=0}^{n-1} (1+t^{2j}) = \frac{2}{1+t^{2n}} \, \prod_{j=1}^{n} (1+t^{2j}), \\ \\
    (-t^{-(2n-1)};t^2)_n &=& \displaystyle \prod_{j=0}^{n-1} (1+t^{2(j-n)+1}) = \displaystyle \prod_{j=0}^{n-1} (1+t^{2(n-j)-1}) \, t^{-2(n-j)+1} \\ 
    &=& \displaystyle \prod_{j=0}^{n-1} (1+t^{2j+1})  \, \prod_{j=0}^{n-1} t^{-2(n-j)+1}, \\ \\
    (-t;t^2)_n &=& \displaystyle \prod_{j=0}^{n-1} (1+t^{2j+1}), \\ \\
    (-t^{-2n};t^2)_n &=& \displaystyle \prod_{j=0}^{n-1} (1+t^{2(j-n)}) = \displaystyle \prod_{j=0}^{n-1} (1+t^{2(n-j)}) \, t^{-2(n-j)} \\ 
    &=& \displaystyle  \prod_{j=1}^{n} (1+t^{2j}) \; \prod_{j=0}^{n-1} t^{-2(n-j)},
\end{array}
$$
we reach the simplified expression
$$
{}_3\phi_2 \left[ 
\begin{array}{c}
t^{-2n}, \ t^{2n}, \ t\\
-t, \ -t^2  
\end{array}
; t^2, t^2
\right] 
\;=\;
\displaystyle \frac{2}{1+t^{2n}} \; \prod_{j=0}^{n-1}t 
\;=\;
\frac{2 \; t^n}{1+t^{2n}} 
$$
that coincides with the claimed result.
\end{proof}

\medskip
\begin{rem}\label{cor:recTn}
In light of the fact that big q-Jacobi polynomials $P_n(z;a,b,c;q)$ are defined as in \eqref{eq:bigq-Jac}  (see also \cite[Equation (7.3.10)]{GasperRahman} and \cite[Equation (1.1)]{GIM92}), 
Prop. \ref{prop:1_su_Tn} leads to the novel identity
$\big(T_n(x) \big)^{-1}=P_n\left(t;-t^{-1}, -t^{-1}, -1; t^2\right)$
for $x$ in \eqref{eq:xvst}. Moreover, comparing Prop. \ref{prop:1_su_Tn} with equation  \eqref{eq:Tn_vs_2F1}, we are also able to establish that
$$
{}_3\phi_2 \left[ 
\begin{array}{c}
t^{-2n}, \ t^{2n}, \ t\\
-t, \ -t^2  
\end{array}
; t^2, t^2
\right] \; {}_2F_1 \Big(-n, n; \frac{1}{2}; -\frac{(t-1)^2}{4t}  \Big)=1.
$$ 
\end{rem}

\medskip
Next we prove that the q-hypergeometric function ${}_3\phi_2$ of Prop. \ref{prop:1_su_Tn} satisfies also another identity.

\medskip
\begin{lem}\label{prop:3Phi2_sum}
For all $n \in \mathbb{N}$,
$$
\begin{array}{l}
(t^{2n}-1)^2 \; \sum_{i=1}^{n} (-1)^i \, \frac{\displaystyle \prod_{\ell=1}^{i} (t^{2\ell-1}-1) \; \prod_{\ell=1}^{i-1} (t^{2(n+\ell)}-1) \;  \prod_{\ell=1}^{i-1} (t^{2(n-\ell)}-1)}{
t^{(2n-i)i-i} \; \displaystyle \prod_{\ell=1}^{i} (t^{2\ell}-1) \; \prod_{\ell=1}^{2i} (t^{\ell}+1)}\;=\\ \\
\qquad =\; {}_3\phi_2 \left[ 
\begin{array}{c}
t^{-2n}, \ t^{2n}, \ t\\
-t, \ -t^2  
\end{array}
; t^2, t^2
\right]-1.
\end{array}
$$
\end{lem}

\begin{proof}
We start expressing the terms in the summation in terms of q-Pochhammer symbols as
\begin{itemize}

    \item[\phantom{1}] 
        $\begin{array}{rcl}
            \displaystyle \prod_{\ell=1}^{i} (t^{2\ell-1}-1) 
            &=& 
            (-1)^i \; \prod_{j=0}^{i-1} \big(1-t \, (t^2)^j\big)
            \;=\; (-1)^i \; (t;t^2)_i, \\ \\
        \end{array}$ 

    \item[\phantom{2}]
        $\begin{array}{rcl}
            \displaystyle \prod_{\ell=1}^{i-1} (t^{2(n+\ell)}-1) &=& (-1)^{i-1} \, \prod_{\ell=1}^{i-1} (1-t^{2(n+\ell)}) = \frac{(-1)^{i-1}}{1-t^{2n}} \, \prod_{\ell=0}^{i-1} \big(1-t^{2n} \, (t^{2})^{\ell}\big) \\
            &=& \frac{(-1)^{i}}{t^{2n}-1} \, (t^{2n};t^2)_i, \\ \\
        \end{array}$ 

    \item[\phantom{3}]
        $\begin{array}{rcl}
            \prod_{\ell=1}^{i-1} (t^{2(n-\ell)}-1) &=&
            \prod_{\ell=1}^{i-1} (1-t^{-2(n-\ell)}) \; 
            t^{2(n-\ell)} = \displaystyle \frac{\displaystyle \prod_{\ell=0}^{i-1} (1-t^{-2n} (t^2)^{\ell})}{1-t^{-2n}} \;  \; \frac{t^{2n(i-1)}}{\displaystyle \prod_{\ell=1}^{i-1} t^{2\ell}} \\
            &=& \frac{(t^{-2n};t^2)_i \; t^{2n(i-1)-i(i-1)}}{1-t^{-2n}},\\ \\
        \end{array}$ 

    \item[\phantom{4}]
        $\begin{array}{rcl}
            \prod_{\ell=1}^{i} (t^{2\ell}-1) &=& (-1)^i \; \prod_{\ell=1}^{i} (1-t^{2\ell}) = (-1)^i \; \prod_{\ell=0}^{i-1} \big(1-t^{2} (t^2)^{\ell}\big) \\
            &=& (-1)^i \; (t^2;t^2)_i, \\ \\
        \end{array}$ 

    \item[\phantom{5}] 
        $\begin{array}{rcl}
            \prod_{\ell=1}^{2i} (t^{\ell}+1) &=& \prod_{\ell=1}^{i} (1+t^{2\ell-1}) \; 
            (1+t^{2\ell}) = \prod_{\ell=0}^{i-1} \big(1+t (t^2)^{\ell}\big) \;  \prod_{\ell=0}^{i-1} \big(1+t^2 (t^2)^{\ell}\big) \\
            &=& (-t;t^2)_i \; (-t^2;t^2)_i.
        \end{array}$ 
        
\end{itemize}
Now, substituting in the summation, we obtain
$$
\begin{array}{c}
    (t^{2n}-1)^2 \; \sum_{i=1}^{n} (-1)^i \, \frac{\displaystyle \prod_{\ell=1}^{i} (t^{2\ell-1}-1) \; \prod_{\ell=1}^{i-1} (t^{2(n+\ell)}-1) \;  \prod_{\ell=1}^{i-1} (t^{2(n-\ell)}-1)}{
    t^{(2n-i)i-i} \; \displaystyle \prod_{\ell=1}^{i} (t^{2\ell}-1) \; \prod_{\ell=1}^{2i} (t^{\ell}+1)} \;=\quad \\ \\
    \quad=\;\sum_{i=1}^{n} \frac{(t;t^2)_i \; (t^{2n};t^2)_i \; (t^{-2n};t^2)_i}{(-t;t^2)_i \; (-t^2;t^2)_i}
    \; \frac{t^{2i}}{(t^2;t^2)_i}, 
\end{array}
$$
and since 
$$
(t^{-2n};t^2)_i=\prod_{j=0}^{i-1} \big(1-t^{-2n} \, (t^2)^j\big) =0 \qquad \text{for} \qquad i>n,
$$
the last sum is equivalent to the series
$$
\begin{array}{c}
\sum_{i=1}^{\infty} \frac{(t;t^2)_i \; (t^{2n};t^2)_i \; (t^{-2n};t^2)_i}{(-t;t^2)_i \; (-t^2;t^2)_i}
\; \frac{t^{2i}}{(t^2;t^2)_i} \;=\; \qquad \qquad \qquad \qquad \qquad 
\\ \\
 \qquad \qquad \qquad \qquad \qquad \;=\;\left( \sum_{i=0}^{\infty} \frac{(t;t^2)_i \; (t^{2n};t^2)_i \; (t^{-2n};t^2)_i}{(-t;t^2)_i \; (-t^2;t^2)_i}
\; \frac{t^{2i}}{(t^2;t^2)_i} \right) -1
\end{array}
$$
which is exactly  
$$
{}_3\phi_2 \left[ 
\begin{array}{c}
t^{-2n}, \ t^{2n}, \ t\\
-t, \ -t^2  
\end{array}
; t^2, t^2
\right]-1,
$$
so concluding the proof.
\end{proof}

\begin{prop} \label{prop:identity_Tn} 
For all $n\in \mathbb{N}$, 
\begin{equation}\label{eq:key_equation_Tn}
\sum_{i=1}^n \frac{2^i}{T_i(x)+1} \;  \prod_{\ell=0}^{i-1}  C_{\ell,i}(x) \, \frac{\Big(T_{\ell}(x)\Big)^2 - \Big(T_{n}(x)\Big)^2}{\Big( T_{\ell}(x) + 1 \Big)^2}  \;=\; \frac{1-T_n(x)}{2T_n(x)}
\end{equation}   
with 
\begin{equation}\label{eq:Cdef}
C_{\ell,i}(x) \; = \; \frac{ \big( T_{\ell}(x) - T_{\ell+1}(x) \big) \ \big( T_{\ell}(x)+1 \big)}{\big( T_{i}(x) - T_{i+1}(x) \big) \, - \, \big( T_{\ell}(x) - T_{\ell+1}(x) \big)}, \quad  
\begin{array}{l}
    i=1,\ldots,n,\\
    \ell=0,\ldots,i-1.
\end{array}
\end{equation}
\end{prop}

\begin{proof}
By exploiting the notation introduced in \eqref{eq:xvst}, we start writing
$$
T_{\ell}(x) + 1\;=\;\frac{t^{\ell} + t^{-\ell} + 2}{2}\;=\; \frac{(t^{\ell} + 1)^2}{2 t^{\ell}}
$$
and
$$
\begin{array}{c}
T_{\ell}(x) - T_{\ell+1}(x) \;=\; 
\frac{t^\ell+t^{-\ell}}{2}\;-\;\frac{t^{\ell+1}+t^{-(\ell+1)}}{2}
\;=\; -\frac{(t^{2\ell+1}-1)(t-1)}{2 t^{\ell+1}}.
\end{array}
$$
Thus, for $i>\ell$,
$$
\begin{array}{c}
\big( T_{i}(x) - T_{i+1}(x) \big) \, - \, \big( T_{\ell}(x) - T_{\ell+1}(x) \big)
\;=\; -\frac{(t^{i+\ell+1}+1) \; (t^{i-\ell}-1) \; (t-1)}{2 t^{i+1}}.
\end{array}
$$
Now, substituting these expressions into $C_{\ell,i}(x)/\big( T_{\ell}(x) + 1 \big)^2$ we get
$$
\begin{array}{c}
\frac{C_{\ell,i}(x)}{\big( T_{\ell}(x) + 1 \big)^2} \;=\; \frac{2 t^i (t^{2\ell+1}-1)}{(t^{i+\ell+1}+1) \; (t^{i-\ell}-1) \; (t^{\ell}+1)^2}.
\end{array}
$$
Hence
$$
\begin{array}{rcl}
\frac{2^i}{T_i(x)+1} \, \prod_{\ell=0}^{i-1} \frac{C_{\ell,i}(x)}{\big( T_{\ell}(x) + 1 \big)^2} 
&=&
\frac{2^{i+1} t^{i}}{(t^i+1)^2} \, \prod_{\ell=0}^{i-1}   \frac{2t^i \, (t^{2\ell+1}-1)}{(t^{i+\ell+1}+1) \; (t^{i-\ell}-1) \; (t^{\ell}+1)^2} \\ \\
&=& \frac{2^{2i+1} t^{i(i+1)}  \; \displaystyle \prod_{\ell=0}^{i-1} (t^{2\ell+1}-1) }{\displaystyle \prod_{\ell=i+1}^{2i} (t^{\ell}+1) \;  \displaystyle \prod_{\ell=1}^{i} (t^{\ell}-1) \; \displaystyle \prod_{\ell=0}^{i} (t^{\ell}+1)^2 } \\  \\
&=& \frac{2^{2i+1} t^{i(i+1)}  \; \displaystyle \prod_{\ell=0}^{i-1} (t^{2\ell+1}-1) }{\displaystyle \prod_{\ell=0}^{2i} (t^{\ell}+1) \; \displaystyle \prod_{\ell=1}^{i} (t^{\ell}-1) \; \displaystyle \prod_{\ell=0}^{i} (t^{\ell}+1) }\\ \\
&=& \frac{2^{2i-1} t^{i(i+1)}   \; \displaystyle \prod_{\ell=1}^{i} (t^{2\ell-1}-1) }{ \displaystyle \prod_{\ell=1}^{2i} (t^{\ell}+1) \; \displaystyle \prod_{\ell=1}^{i} (t^{\ell}-1) \; \displaystyle \prod_{\ell=1}^{i} (t^{\ell}+1) }  \\ \\
&=& \frac{2^{2i-1}   \; \displaystyle \prod_{\ell=1}^{i} (t^{2\ell-1}-1) }{t^{-i(i+1)} \; \displaystyle \prod_{\ell=1}^{2i} (t^{\ell}+1) \; \displaystyle \prod_{\ell=1}^{i} (t^{2\ell}-1) }. 
\end{array}
$$
Moreover, for $n>\ell$,
$$
\begin{array}{rcl}
\Big(T_{\ell}(x)\Big)^2 - \Big(T_{n}(x)\Big)^2 
&=& -\frac{(t^{2(n+\ell)}-1) (t^{2(n-\ell)}-1)}{2^2 t^{2n}}
\end{array}
$$
and hence
$$
\begin{array}{rcl}
\prod_{\ell=0}^{i-1} \Big( \big(T_{\ell}(x)\big)^2 - \big(T_{n}(x)\big)^2 \Big)
&=& \frac{(-1)^i}{2^{2i} t^{2ni}} \; \prod_{\ell=0}^{i-1} (t^{2(n+\ell)}-1) (t^{2(n-\ell)}-1) \\ \\
&=&  \frac{(-1)^i \; (t^{2n}-1)^2}{2^{2i} t^{2ni}} \; \prod_{\ell=1}^{i-1} (t^{2(n+\ell)}-1) (t^{2(n-\ell)}-1).
\end{array}
$$
Then, it follows that
$$
\begin{array}{l}
\sum_{i=1}^n \frac{2^i}{T_i(x)+1} \;  \prod_{\ell=0}^{i-1}  C_{\ell,i}(x) \, \frac{\Big(T_{\ell}(x)\Big)^2 - \Big(T_{n}(x)\Big)^2}{\Big( T_{\ell}(x) + 1 \Big)^2}
\;= \\ \\
\quad=\;(t^{2n}-1)^2 \; \sum_{i=1}^n (-1)^i \; 
\frac{2^{2i-1}  \; \displaystyle \prod_{\ell=1}^{i} (t^{2\ell-1}-1) \; 
\displaystyle 
\prod_{\ell=1}^{i-1} (t^{2(n+\ell)}-1)  \;  
(t^{2(n-\ell)}-1) }
{2^{2i} \; t^{2ni-i(i+1)} \; \displaystyle \prod_{\ell=1}^{2i} (t^{\ell}+1) \; \displaystyle \prod_{\ell=1}^{i} (t^{2\ell}-1) }  
\\ \\
\quad=\;\frac{1}{2} \; (t^{2n}-1)^2 \; \sum_{i=1}^{n} (-1)^i \, \frac{\displaystyle \prod_{\ell=1}^{i} (t^{2\ell-1}-1) \; \prod_{\ell=1}^{i-1} (t^{2(n+\ell)}-1) \;  \prod_{\ell=1}^{i-1} (t^{2(n-\ell)}-1)}{
t^{(2n-i)i-i} \; \displaystyle \prod_{\ell=1}^{i} (t^{2\ell}-1) \; \prod_{\ell=1}^{2i} (t^{\ell}+1)}. 
\end{array}
$$
Since, in light of Proposition \ref{prop:1_su_Tn} and Lemma \ref{prop:3Phi2_sum} we can write 
$$
\begin{array}{c}
(t^{2n}-1)^2 \; \sum_{i=1}^{n} (-1)^i \, \frac{\displaystyle \prod_{\ell=1}^{i} (t^{2\ell-1}-1) \; \prod_{\ell=1}^{i-1} (t^{2(n+\ell)}-1) \;  \prod_{\ell=1}^{i-1} (t^{2(n-\ell)}-1)}{
t^{(2n-i)i-i} \; \displaystyle \prod_{\ell=1}^{i} (t^{2\ell}-1) \; \prod_{\ell=1}^{2i} (t^{\ell}+1)} \;=\qquad \\ \\
\qquad =\; \frac{2 t^n}{1+t^{2n}} -1 \;=\; -\frac{(1-t^n)^2}{1+t^{2n}}\;=\;\frac{1-T_n(x)}{T_n(x)},
\end{array}
$$
we conclude that
$$
\sum_{i=1}^n \frac{2^i}{T_i(x)+1} \; \prod_{\ell=0}^{i-1}  C_{\ell,i}(x) \, \frac{\Big(T_{\ell}(x)\Big)^2 - \Big(T_{n}(x)\Big)^2}{\Big( T_{\ell}(x) + 1 \Big)^2} \;=\; \frac{1-T_n(x)}{2T_n(x)}.
$$
\end{proof}

\section{State-of-the-art on minimum-support interpolating subdivision schemes
reproducing finite sets of integer powers of exponentials}\label{sec:0}

Subdivision schemes are efficient iterative algorithms capable of rapidly and accurately visualizing planar and spatial curves starting from a coarse set of points ${\bf p}^{(0)}$  enclosing the shape of the desired curve. 
In a nutshell, for $k\in \mathbb{N}_0$, the iterative algorithm computes the set of points ${\bf p}^{(k+1)}$ by means of local linear combinations of the set of points ${\bf p}^{(k)}$ with coefficients given by the entries of the so-called $k$-level \emph{subdivision mask}. Very frequently the $k$-level subdivision mask is replaced with the Laurent polynomial defined by its entries, which we formally call the $k$-level \emph{subdivision symbol}. The following theorem is a special case of the general result introduced in \cite{CR2011} to characterize (in terms of the algebraic properties of the $k$-level symbol and its first derivative) the capability of an odd-symmetric, convergent, non-singular (i.e., converging to the constant zero function if and only if applied to the initial sequence of all zero values) subdivision scheme of generating/reproducing the function space 
$$
V_{2n+2,\theta}\;=\;{\rm span}\left\{\;1,\;(\bigcdot), \;  \{e^{\pm \mathrm{i} j\theta \, \bigcdot}\}_{j=1}^{n}\;\right\}\;\cap\;\mathcal{C}^\infty(\mathbb{R},\mathbb{R}),$$
with $n \in \mathbb{N}, \ \theta \in [0, {\sigma_{\star}(n)}[\; \cup\; \mathrm{i} \mathbb{R}^+$
and $\sigma_{\star}(n)\;\in\;]0,2\pi[$ the so-called critical length of the space (see, e.g., \cite{BCM19,BCM20,CMP04}). 
{ 
The critical length is the number identifying
the maximum width of the parameter $\theta$ domain for which any associated spline space admits a B-spline basis.
}
With regards to the family of exponential B-spline schemes generating $V_{2n+2,\theta}$, the $k$-level symbol has a well-known expression (see, e.g., \cite{CGR16,Yoon13,Jeong13,Romani09}) that can be concisely written as
\begin{equation}\label{eq:sigma_symbol_via_alfa}
s_{2n+2,k}(z) \;=\;
2 \prod_{\ell=0}^{n} \ a_{\ell,k}(z)
\end{equation}
with
\begin{equation} \label{eq:az}
    a_{\ell,k}(z)\;=\;\frac{z + 2 \, T_{\ell}\big( v_{k} \big) +z^{-1}}{2 \Big( T_{\ell}\big( v_{k} \big)+1 \Big)}, \qquad \ell=0,\ldots,n
\end{equation}
and  
\begin{equation}\label{rec_par}
v_{k}\;=\;\cos\left( \frac{\theta}{2^{k+1}} \right).
\end{equation} 

\begin{rem}\label{rem_bspline_similarity}
We observe that $a_{0,k}(z)\;=\;\frac{(1+z)^{2}}{4z}$. Moreover, $\lim_{k\rightarrow +\infty}
v_k=1$ and, for all $\ell\in\mathbb{N}$,
$$
\lim_{k\rightarrow +\infty}
T_{\ell}(v_k)\;=\;1, \qquad  
\lim_{k\rightarrow +\infty} a_{\ell,k}(z)\;=\;\frac{(1+z)^{2}}{4z}, 
$$
so that 
$$\lim_{k\rightarrow +\infty} s_{2n+2,k}(z)\;=\;\frac{(1+z)^{2n+2}}{2^{2n+1} \, z^{n+1}},$$ 
which means that, as $k \rightarrow +\infty$, the exponential B-spline scheme with $k$-level symbol $s_{2n+2,k}(z)$ converges to the polynomial B-spline scheme of degree $2n+1$ (see also \cite[Prop. 6.4]{CGR16}).
\end{rem}

\begin{thm}{\cite[Theorem 1]{CR2011}}\label{theo_recalls} 
The non-singular subdivision scheme with odd-symmetric symbols $\{a_k(z)\}_{k \in \mathbb{N}_0}$ 
(i.e. satisfying $a_k(z)=a_k(z^{-1})$ for all $k \in \mathbb{N}_0$)
generates functions in $V_{2n+2,\theta}$  if and only if, for each $k \in \mathbb{N}_0$ and for $r_k=e^{\mathrm{i} \frac{\theta}{2^{k+1}}}$, the $2n+2$ conditions 
\begin{equation}\label{eq:Vgen_condiz}
\begin{array}{l}
a_k(-1)=0, \\
a_k'(-1)=0,\\
a_k(-r_k^{\pm j})=0, \quad \hbox{for all} \quad j=1,...,n  
\end{array}
\end{equation}
are all satisfied.
Furthermore, the subdivision scheme reproduces functions in  $V_{2n+2,\theta}$ if and only if, for each $k \in \mathbb{N}_0$, the $2n+2$ conditions 
\begin{equation}\label{eq:Vrep_condiz}
\begin{array}{l}
a_k(1)=2, \\
a_k'(1)=0,\\
a_k(r_k^{\pm j})=2, \quad \hbox{for all} \quad j=1,...,n 
\end{array}
\end{equation}
are also all satisfied.    
\end{thm}

\begin{rem}
The first two conditions in \eqref{eq:Vgen_condiz} are for the generation of $1,(\bigcdot)$ and the last  $2n$ for the generation of $\{e^{\pm \mathrm{i} j\theta \, \bigcdot}\}_{j=1}^{n}$. An analogous observation holds for \eqref{eq:Vrep_condiz} and the reproduction property.
\end{rem}

\smallskip
We additionally recall that a non-stationary subdivision scheme 
with odd-symmetric symbols $\{a_k(z)\}_{k \in \mathbb{N}_0}$ is termed \emph{interpolatory} if and only if, for all $k \in \mathbb{N}_0$, its $k$-level symbol satisfies the property (see, e.g., \cite[eq.(2.7)]{CGR11})  
\begin{equation}\label{eq:interp_condiz}
a_k(z)+a_k(-z)\;=\;2.
\end{equation}
In light of Theorem \ref{theo_recalls}, it is easy to see that, if an interpolatory scheme generates $V_{2n+2,\theta}$, then it also reproduces it. Indeed, if \eqref{eq:interp_condiz} and \eqref{eq:Vgen_condiz} hold true, then all conditions in \eqref{eq:Vrep_condiz} are fulfilled too. The class of minimum-support interpolating subdivision schemes
reproducing $V_{2n+2,\theta}$
was investigated in several previous papers (see, e.g., \cite{CCR,CGR11,CGR16,CR2010,DLL03,NR15,Romani09}) proving their convergence and smoothness properties, but, for arbitrary values of $n$ their $k$-level symbol (hereinafter denoted by $m_{2n+2,k}(z)$) is still not known in explicit closed-form and its expression has to be determined by means of a suitable algorithm aimed at calculating the Laurent polynomial (say $\ell_k(z)$) to be multiplied by the exponential B-spline symbol $s_{2n+2,k}(z)$ generating the same function space \cite{CGR11}. More in detail, the computational procedure currently used to build $m_{2n+2,k}(z)$ reads as follows:

\begin{itemize}

    \item[(i)] let $s_{2n+2,k}(z)$ be the $k$-level symbol of the exponential B-spline scheme generating $V_{2n+2,\theta}$;

    \item[(ii)] construct the $(2n+1) \times (2n+1)$ matrix ${\bf A}_k$, leading principal submatrix of a Hurwitz type matrix associated to $s_{2n+2,k}(z)$;

    \item[(iii)] determine $({\bf A}_k)^{-1}$;

    \item[(iv)] set $\boldsymbol{\ell}_k$ to be the $(n+1)$st row of $({\bf A}_k)^{-1}$ and build the Laurent polynomial $\ell_k(z)=\sum_{j=-n}^{n} (\boldsymbol{\ell}_{k})_j \, z^j$;

    \item[(v)] construct the interpolatory symbol of the minimum-support interpolating scheme reproducing $V_{2n+2,\theta}$ as $m_{2n+2,k}(z)=s_{2n+2,k}(z) \ \ell_{k}(z)$.

\end{itemize}

\begin{exmp}[$n=3$] 
The $k$-level symbol of the exponential B-spline scheme generating $ 
V_{8,\theta}=\hbox{span} \big \{ \ 1, (\bigcdot),  e^{\pm {\rm i} \theta \bigcdot}, e^{\pm 2{\rm i} \theta \bigcdot}, e^{\pm 3{\rm i}\theta \bigcdot}   \ \big \}
$ is $s_{8,k}(z)= \sum_{j=0}^8 s_j(v_k) \, z^{j-4}$ with
$$
\begin{array}{c}
s_0(v_k)\;=\;s_8(v_k)\;=\; \frac{1}{32 v_k^2 (2v_k - 1)^2 (v_k + 1)^2},\\ \\
s_1(v_k)\;=\;s_7(v_k)\;=\; \frac{1}{8 v_k (2v_k^2 + v_k - 1)}, \\ \\
s_2(v_k)\;=\;s_6(v_k)\;=\; \frac{8 v_k^3 + 4 v_k^2 - 4v_k - 1}{8(2v_k^2 + v_k - 1)^2},\\ \\
s_3(v_k)\;=\;s_5(v_k)\;=\; \frac{8 v_k^3 + 4 v_k^2 - 4v_k - 1}{8 v_k (2v_k^2 + v_k - 1)}, \\ \\
s_4(v_k)\;=\; \frac{32 v_k^5 + 32 v_k^4 - 16 v_k^3 - 16 v_k^2 + 2 v_k + 1}{16 v_k^2 (2v_k - 1) (v_k + 1)^2}.
\end{array}
$$
The matrix ${\bf A}_k$ extracted from the Hurwitz type matrix associated to $s_{8,k}(z)$ is
$$
{\bf A}_k=\left [ 
\begin{matrix}
s_1(v_k) & s_3(v_k) & s_5(v_k) & s_7(v_k) & 0        & 0        & 0 \cr 
s_0(v_k) & s_2(v_k) & s_4(v_k) & s_6(v_k) & s_8(v_k) & 0        & 0 \cr 
0        & s_1(v_k) & s_3(v_k) & s_5(v_k) & s_7(v_k) & 0        & 0 \cr 
0        & s_0(v_k) & s_2(v_k) & s_4(v_k) & s_6(v_k) & s_8(v_k) & 0 \cr 
0        & 0        & s_1(v_k) & s_3(v_k) & s_5(v_k) & s_7(v_k) & 0 \cr 
0        & 0        & s_0(v_k) & s_2(v_k) & s_4(v_k) & s_6(v_k) & s_8(v_k) \cr 
0        & 0        & 0        & s_1(v_k) & s_3(v_k) & s_5(v_k) & s_7(v_k) 
\end{matrix}
\right ],
$$
and the coefficient vector $\boldsymbol{\ell}_{k}$ for constructing the Laurent polynomial $\ell_{k}(z)$ and hence the interpolatory symbol  
$m_{8,k}(z)=s_{8,k}(z) \ell_{k}(z)$ is 
the 4th row of $({\bf A}_k)^{-1}$.  
\end{exmp}

\section{Closed-form symbols of minimum-support interpolating subdivision schemes reproducing finite sets of integer powers of exponentials}\label{sec:3}

This section is devoted to present our novel results in the domain of subdivision schemes. 
They consist in the derivation of a closed-form expression (that avoids using the algorithm recalled in Section \ref{sec:0}) for the
$k$-level symbol of the minimum-support  interpolating subdivision scheme reproducing $V_{2n+2,\theta}$ 
(Theorem \ref{theo_main}) and in showing that this scheme is the non-stationary version of the $(2n+2)$-point Dubuc-Deslauriers interpolating scheme reproducing algebraic polynomials whose closed-form expression was firstly presented in \cite{DS2006} (Proposition \ref{prop_limit}).  

{ 
\begin{thm}\label{theo_main} 
Let $n\in\mathbb{N}_0$ and $V_{2n+2,\theta}$ as in \eqref{eq:V_space}. Define $v_k=\cos(2^{-k-1}\theta)$,  
\begin{equation}\label{eq:small_civk}
    c_i(v_k)\;=\;\frac{2^i}{T_i(v_{k}) +1} \, \displaystyle \prod_{\ell=0}^{i-1} C_{\ell,i}(v_{k}),\qquad i=1, \ldots, n,
\end{equation}
and 
\begin{equation}\label{eq:bz_def} 
    b_{i,k}(z)\;=\;\big(2 a_{0,k}(z)-1\big) \, c_i(v_k)  \, \prod_{\ell=0}^{i-1} a_{\ell,k}(-z) 
\end{equation}  
with $C_{\ell,i}(v_{k})$ and $a_{\ell,k}(z)$ as in \eqref{eq:Cdef} and \eqref{eq:az}, respectively. 
Then, the Laurent polynomial 
\begin{equation}\label{eq:a2np2}
m_{2n+2,k}(z) \; = \; 2 a_{0,k}(z) \, \left(   1+ \sum_{i=1}^n b_{i,k}(z) \; \prod_{\ell=1}^{i-1} a_{\ell,k}(z)  \right)   
\end{equation}   
is the $k$-level symbol of the $(2n+2)$-point non-stationary odd-symmetric interpolatory subdivision scheme reproducing $V_{2n+2,\theta}$.
\end{thm}

\begin{proof} \,
We start observing that, according to \eqref{eq:a2np2}, 
{ 
 for all $k \in \mathbb{N}_0$,
$z^{2n+1} \, m_{2n+2,k}(z)$ is a polynomial of degree $2(2n+1)$, namely the non-zero coefficients of
the Laurent polynomial $m_{2n+2,k}(z)$ are  
indexed between $-(2n+1)$ and $(2n+1)$.
Thus $\{m_{2n+2,k}(z)\}_{k \in \mathbb{N}_0}$ identifies a $(2n+2)$-point subdivision scheme.} Now we continue by proving (a) odd-symmetry, (b) interpolation and, by using induction, (c) generation of $V_{2n+2,\theta}$.
\begin{description} 
        \item[\sl (a) Odd-symmetry:]  
from \eqref{eq:az} we have that $a_{j,k}(z^{-1})=a_{j,k}(z)$, for every $j\in\mathbb{N}_0$. Thus, due to \eqref{eq:bz_def}, we have $b_{j,k}(z^{-1})=b_{j,k}(z)$, for every $j\in\mathbb{N}$. Therefore, from \eqref{eq:a2np2} the $k$-level symbol $m_{2n+2,k}(z)$ satisfies $m_{2n+2,k}(z^{-1})=m_{2n+2,k}(z)$ for all $n\in\mathbb{N}_0$, i.e., 
it defines an odd-symmetric subdivision scheme.\\

        \item[\sl (b) Interpolation:] 
recalling that, for all $k\in\mathbb{N}_0$,
$$
s_{2n,k}(z) \;=\; 2 \prod_{i=0}^{n-1} \ a_{i,k}(z),
$$
we have, according to \eqref{eq:a2np2},
\begin{equation}\label{eq:rel_m2np2_m2n}
    \begin{array}{rcl}
        m_{2n+2,k}(z)\;-\;m_{2n,k}(z) &=& 2 \, a_{0,k}(z) \, b_{n,k}(z) \,  \prod_{\ell=1}^{n-1} a_{\ell,k}(z)\\ \\
        &=& s_{2n,k}(z)  \ b_{n,k}(z),
    \end{array}
\end{equation}
and analogously
            \[
                m_{2n+2,k}(-z)\;-\;m_{2n,k}(-z)\;=\;s_{2n,k}(-z)\;b_{n,k}(-z).
            \]
Summing these two equations, due to \eqref{eq:bz_def}, we see that  
\begin{equation*} \label{eq:temp_step_interp}
    \begin{array}{l}
        \left( m_{2n+2,k}(z)+ m_{2n+2,k}(-z) \right) - \left(m_{2n,k}(z) + m_{2n,k}(-z)\right)  \;=\\ \\
        \qquad=\; s_{2n,k}(z) \, b_{n,k}(z) \;+\; s_{2n,k}(-z) \, b_{n,k}(-z)\\ 
        \qquad=\;\Big( a_{0,k}(z) + a_{0,k}(-z) - 1\Big) \
        \frac{2^{n+2}}{T_n(v_{k})+1} \ \prod_{i=0}^{n-1}  C_{i,n}(v_{k}) \  a_{i,k}(-z) \, a_{i,k}(z). \\
    \end{array}
\end{equation*}
Then, in light of the fact that 
\[
a_{0,k}(z) + a_{0,k}(-z) -1\;=\;\frac{(1+z)^2}{4z}-\frac{(1-z)^2}{4z} -1  \;=\;0,
\]
we obtain
\[
    m_{2n+2,k}(z)\;+\;m_{2n+2,k}(-z)\;=\;m_{2n,k}(z)\;+\;m_{2n,k}(-z).
\]
Applying the same reasoning recursively, we conclude that
\begin{equation*}\label{eq:chain_eqs}
    \begin{array}{rcl}
    m_{2n+2,k}(z)\;+\;m_{2n+2,k}(-z)&=&m_{2,k}(z)\;+\;m_{2,k}(-z)\\ \\
&\overset{\mathrm{\eqref{eq:a2np2} \ with \ } n=0}{=}&
    2\left(\;a_{0,k}(z)\;+\;a_{0,k}(-z)\;\right)\\ \\
    &=& 
    \frac{(1+z)^2}{2z} - \frac{(1-z)^2}{2z}
    \;=\;2,
    \end{array}
\end{equation*}
which (according to \eqref{eq:interp_condiz}) means that the scheme with $k$-level symbol $m_{2n+2,k}(z)$ in \eqref{eq:a2np2} is interpolating. \\ 

        \item[\sl (c) Generation of $V_{2n+2,\theta}$:]  
            according to Theorem \ref{theo_recalls}, for the generation of the two monomials in $V_{2n+2,\theta}$ we only need to check that, for all $k\in\mathbb{N}_0$,
            $$ 
            m_{2n+2,k}(-1)\;=\;m_{2n+2,k}'(-1)\;=\;0.
            $$
            This simply follows from \eqref{eq:a2np2} and the fact that $a_{0,k}(z)$ has $-1$ as a double root.\\
            In order to guarantee also the generation of $\{e^{\pm \mathrm{i} j\theta \, \bigcdot}\}_{j=1}^{n}$, we further need to prove that, for all $j\in\{1,\dots,n\}$,
            \begin{equation} \label{eq:all_roots}
                m_{2n+2,k}(-r^{\pm j}_{k})\;=\;0, \qquad \text{with} \qquad r_{k}=e^{\mathrm{i} \frac{\theta}{2^{k+1}}}.
            \end{equation}
             To show \eqref{eq:all_roots} we  proceed inductively. \\
        
            {$\bullet$ \, \it Base step $(n=1)$.} We start observing that 
            $\{-r_k^{\pm 1}\}$
            are roots of $m_{4,k}(z)$. Indeed, in light of \eqref{eq:a2np2} with $n=1$, $m_{4,k}(z)$ can be explicitly written as 
            $$m_{4,k}(z)\;=\;2 a_{0,k}(z) \, \left(   1+   b_{1,k}(z)   \right) \;=\; 2 a_{0,k}(z) \, \left(\;1\;-\;\frac{(z^2+1)(z-1)^2}{4z^2\;v_k\;\left(\;1+v_k\right)\;}\;\right),$$
            and, since $(r_k^{\pm 2}+1)/(2r_k^{\pm 1})=v_k$, it follows that
 \[
 \begin{array}{c}
    \frac{(z^2+1)(z-1)^2}{4z^2} \Big \vert_{z=-r_{k}^{\pm 1}} \;=\; \frac{((-r_{k}^{\pm 1})^2+1) (-r_{k}^{\pm 1}-1)^2}{4(-r_{k}^{\pm 1})^2}\qquad \qquad\\ \\
    \qquad \qquad=\;v_k \frac{r_{k}^{\pm 2}+2r_{k}^{\pm 1}+1}{2r_{k}^{\pm 1}} 
    \;=\; v_k \; \left( 1 + \frac{r_k^{\pm 2}+1}{2r_k^{\pm 1}} \right)
    \;=\;v_k\;\left(\;1+v_k\;\right),
 \end{array}
 \]
 thus implying $m_{4,k}(-r_k^{\pm 1})=0$.\\

            $\bullet$ \, {\it Inductive step.} We now suppose (inductive assumption)
            \[
                m_{2(n-1)+2,k}(-r^{\pm j}_{k})\;=\;0,\qquad \forall j\in\{1,\dots,n-1\},
            \]
            and we want to prove that the same holds for $n$ in place of $n-1$. Due to \eqref{eq:rel_m2np2_m2n} and the fact that, for all $j\in\{1,\dots,n-1\}$, $s_{2n,k}(-r^{\pm j}_{k})=0$, we straightforwardly have
            \[
                m_{2n+2,k}(-r^{\pm j}_{k})\;=\;m_{2n,k}(-r^{\pm j}_{k})\;+\;s_{2n,k}(-r^{\pm j}_{k})\;b_{n,k}(-r^{\pm j}_{k})\;=\;0.
            \]
Then, we only need to prove $m_{2n+2,k}(-r^{\pm n}_{k})=0$.   
After observing that
\begin{equation}\label{eq:Tn_rn}
    \frac{r_{k}^{\pm n} + (r_{k}^{\pm n})^{-1}}{2}\;=\; \cos( n \, \theta/2^{k+1})\;=\;T_n(v_k),
\end{equation}
it follows
\begin{equation}\label{eq:prod_aellval}
\begin{array}{c}
a_{\ell,k}(r_{k}^{\pm n}) a_{\ell,k}(-r_{k}^{\pm n}) \;=\;\frac{\Big(T_{\ell}(v_k)\Big)^2 - \Big(T_{n}(v_k)\Big)^2}{\Big( T_{\ell}(v_k) + 1 \Big)^2},
\end{array}
\end{equation}
and so
\begin{equation} \label{eq:relazione_chiave}
\begin{array}{l}
\hspace{-1.0cm} \sum_{i=1}^n c_i(v_k) \;
 \prod_{\ell=0}^{i-1} a_{\ell,k}(r_{k}^{\pm n}) \,  a_{\ell,k}(-r_{k}^{\pm n})\; =\\ \\
\overset{\mathrm{\eqref{eq:small_civk}}}{=}\;
\sum_{i=1}^n \frac{2^i}{T_i(v_k)+1} \; \prod_{\ell=0}^{i-1} C_{\ell,i}(v_{k}) \, a_{\ell,k}(r_{k}^{\pm n}) \,  a_{\ell,k}(-r_{k}^{\pm n}) \\ \\
\overset{\mathrm{\eqref{eq:prod_aellval}}}{=}\;
\sum_{i=1}^n \frac{2^i}{T_i(v_k)+1} \; \prod_{\ell=0}^{i-1} C_{\ell,i}(v_{k}) \, \frac{\Big(T_{\ell}(v_k)\Big)^2 - \Big(T_{n}(v_k)\Big)^2}{\Big( T_{\ell}(v_k) + 1 \Big)^2} \\ \\
\overset{\mathrm{\eqref{eq:key_equation_Tn}}}{=}\;
\frac{1-T_n(v_k)}{2 T_n(v_k)} \;
\overset{\mathrm{\eqref{eq:Tn_rn}}}{=} \;
  -\frac{( 1-r_{k}^{\pm n} )^2}{2 \, (1+r_{k}^{\pm 2n})} \;
  \overset{\mathrm{\eqref{eq:az}}}{=} \;
 \frac{a_{0,k}(-r_{k}^{\pm n})}{1-2a_{0,k}(-r_{k}^{\pm n})}.
\end{array}
\end{equation}
Now the closed-form expression of $m_{2n+2,k}(z)$ in \eqref{eq:a2np2} allows us writing
$$
\begin{array}{l}
m_{2n+2,k}(-r_{k}^{\pm n})\;=\; 2 a_{0,k}(-r_k^{\pm n}) \, \left(   1+ \sum_{i=1}^n b_{i,k}(-r_k^{\pm n}) \; \prod_{\ell=1}^{i-1} a_{\ell,k}(-r_k^{\pm n}) \right)\\
\;\overset{\mathrm{\eqref{eq:bz_def}}}{=}\;
2 a_{0,k}(-r_{k}^{\pm n}) + 
2 \big(2 a_{0,k}(-r_{k}^{\pm n})-1\big) 
\sum_{i=1}^n c_i(v_k) \; \prod_{\ell=0}^{i-1} a_{\ell,k}(r_{k}^{\pm n})  a_{\ell,k}(-r_{k}^{\pm n}) \\
\;\overset{\mathrm{\eqref{eq:relazione_chiave}}}{=}\;
2 a_{0,k}(-r_{k}^{\pm n}) - 
2 a_{0,k}(-r_{k}^{\pm n})=0,
\end{array}
$$ 
thus concluding the proof.
\end{description} 
\end{proof}
}

\begin{rem}
When $n=1,2,3$, the closed-form expression in \eqref{eq:a2np2} gives us back the Laurent polynomials of the known 4-point, 6-point, 8-point interpolatory subdivision schemes proposed in \cite{BCR07,Romani09,CR2010}, respectively.
\end{rem} 

\begin{figure}[h!]
\centering
\includegraphics[scale=0.225]{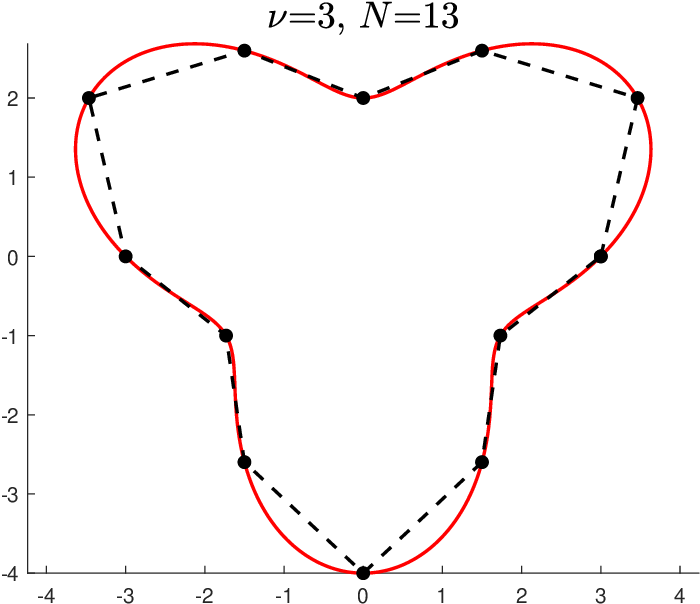}
\includegraphics[scale=0.225]{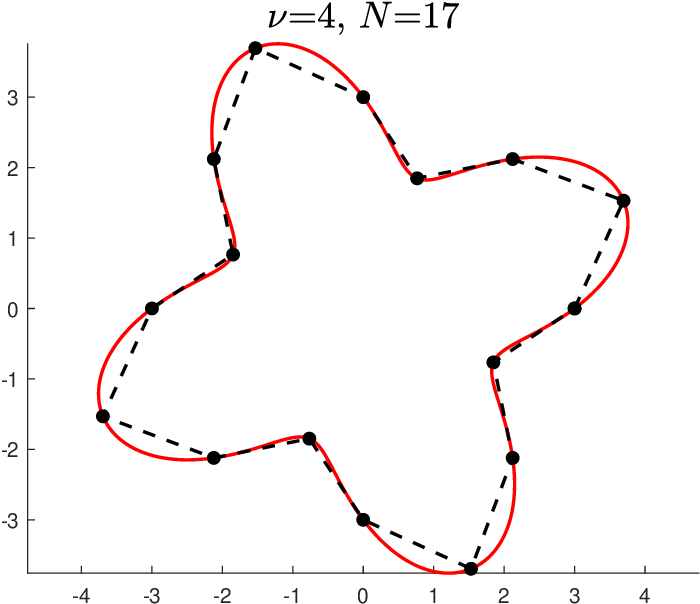}
\includegraphics[scale=0.225]{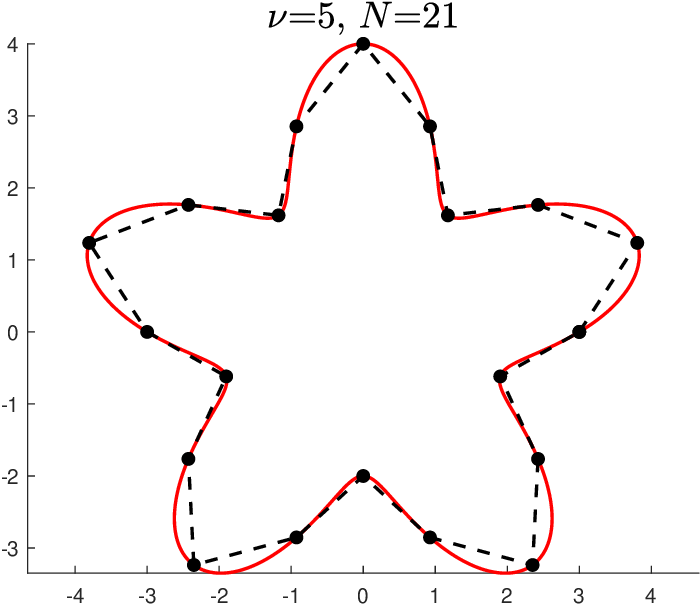}
\includegraphics[scale=0.225]{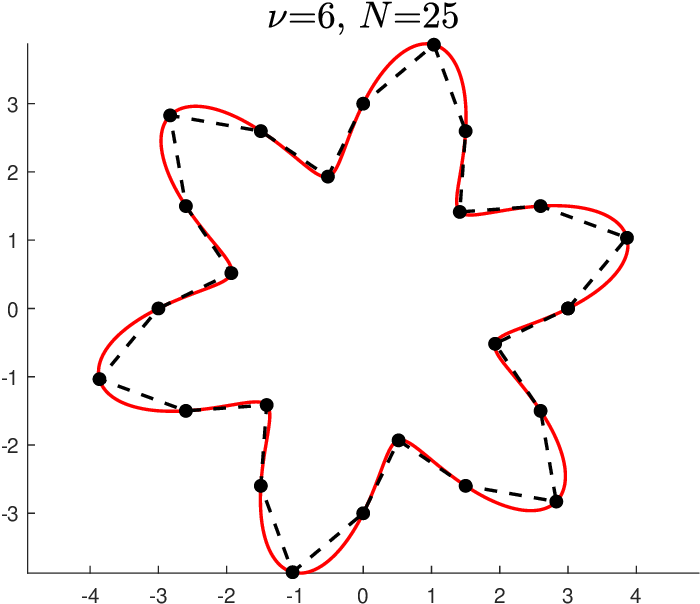}\\
\includegraphics[scale=0.35,trim= 2cm 0 2cm 0,clip]{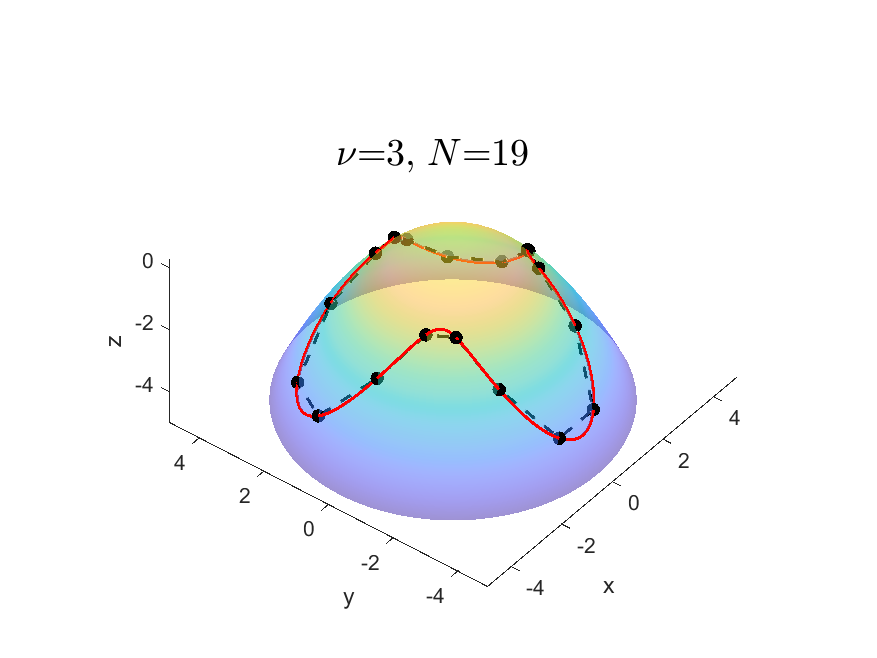} 
\includegraphics[scale=0.35,trim= 2cm 0 2cm 0,clip]{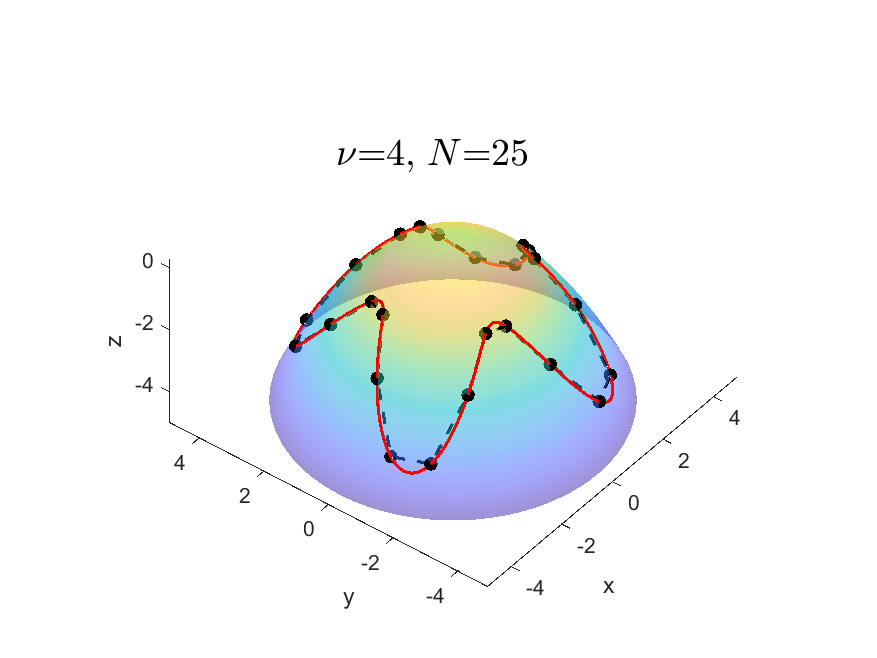} 
\includegraphics[scale=0.35,trim= 2cm 0 2cm 0,clip]{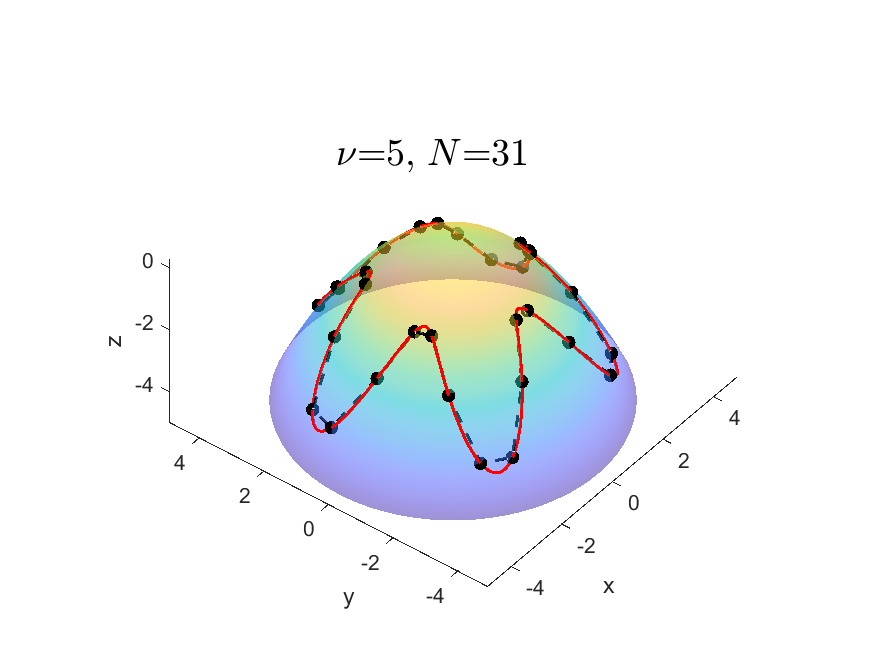} 
\caption{Examples of $\nu$-pointed star-shaped curves
with parametric equations $(x(u),y(u))$ (first row) and $(x(u),y(u),z(u))$ (second row) where
$x(u)=\big( 3+\sin(\nu u) \big) \cos(u)$,  
$y(u)=\big( 3+\sin(\nu u) \big) \sin(u)$, 
$z(u)=-\big(3+\sin(\nu u) \big)^2/4$, $u \in [0, 2\pi]$, obtained after 6 subdivision steps of the interpolatory schemes $m_{2n+2,k}(z)$ reproducing $V_{2n+2,\theta}$ 
with $n=\nu+1$ (first row) while $n=2\nu$ (second row) and $\theta=(2\pi)/(N-1)$,
starting from the dashed polygons with vertices
${\bf p}_i^{(0)}=(x_i^{(0)}, y_i^{(0)})$, $i=1,\ldots,N$ (first row) and 
${\bf p}_i^{(0)}=(x_i^{(0)}, y_i^{(0)}, z_i^{(0)})$, $i=1,\ldots, N$ (second row) where $x_i^{(0)}= \left( 3+\sin \left(\nu \theta (i-1) \right) \right) \cos \left( \theta (i-1) \right)$,  
$y_i^{(0)}=\left( 3+\sin \left( \nu \theta (i-1) \right) \right) \sin \left( \theta (i-1) \right)$,
$z_i^{(0)}=-\left(3+\sin \left(\nu \theta (i-1) \right) \right)^2/4$.}
\label{fig_starshape2d3d}
\end{figure}

\begin{figure}[h!]
\centering
\includegraphics[scale=0.42]{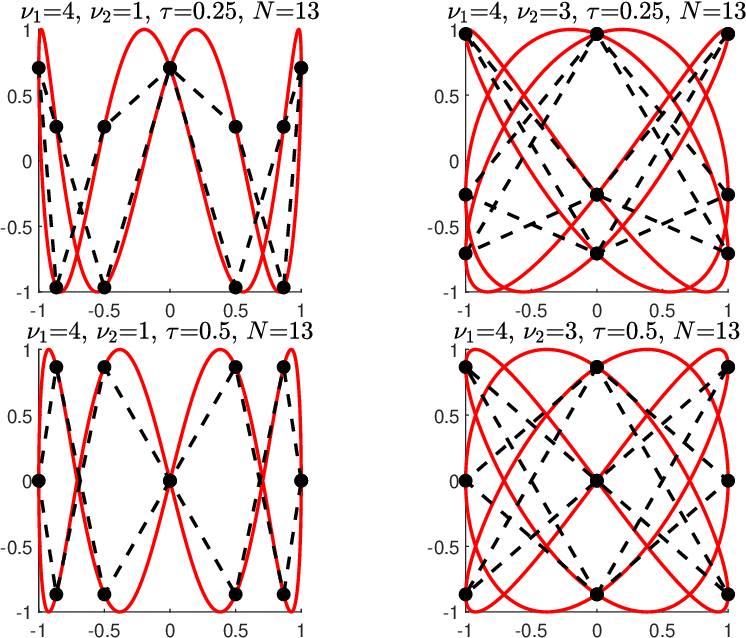} \hspace{0.3cm}
\includegraphics[scale=0.42]{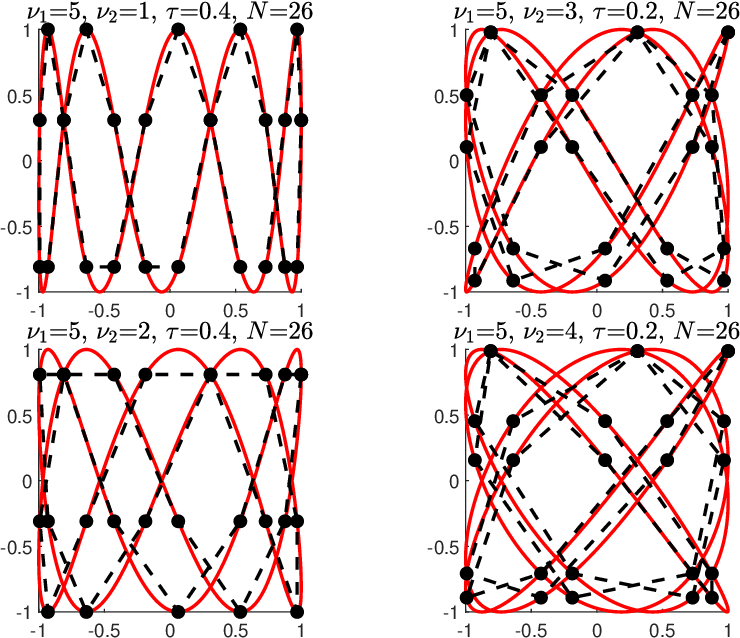}
\caption{Examples of planar Lissajous curves with parametric equations $x(u)= \cos(\nu_2 u)$,  
$y(u)= \cos\big(\nu_1 u-(\tau \pi)/\nu_2 \big)$, with $u \in [0, 2\pi]$, obtained after 6 subdivision steps of the interpolatory schemes $m_{10,k}(z)$ (first two columns) and $m_{12,k}(z)$ (last two columns) reproducing $V_{2n+2,\theta}$ with $\theta=(2\pi)/(N-1)$ and $n=\max\{\nu_1,\nu_2\}$,
starting from the dashed polygons with vertices ${\bf p}_i^{(0)}=(x_i^{(0)}, y_i^{(0)})$, $i=1,\ldots, N$ where $x_i^{(0)}= \cos \left(\nu_2 \theta (i-1) \right)$,  
$y_i^{(0)}= \cos \left(\nu_1 \theta (i-1) -(\tau \pi)/ \nu_2 \right)$.} 
\label{fig_lissajous}
\end{figure}

\begin{figure}[h!]
\centering
\includegraphics[scale=0.28]{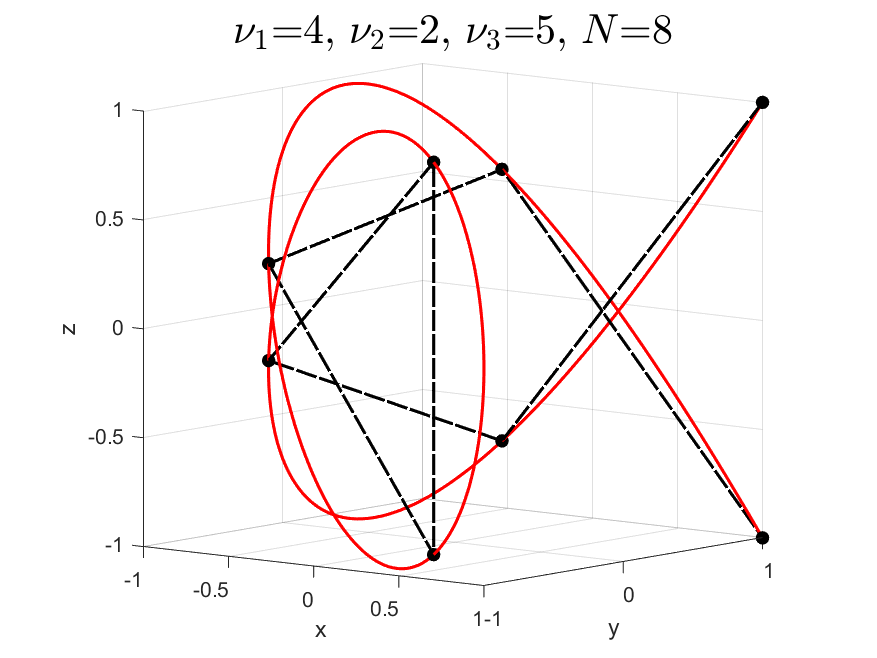}\hspace{-0.5cm}
\includegraphics[scale=0.28]{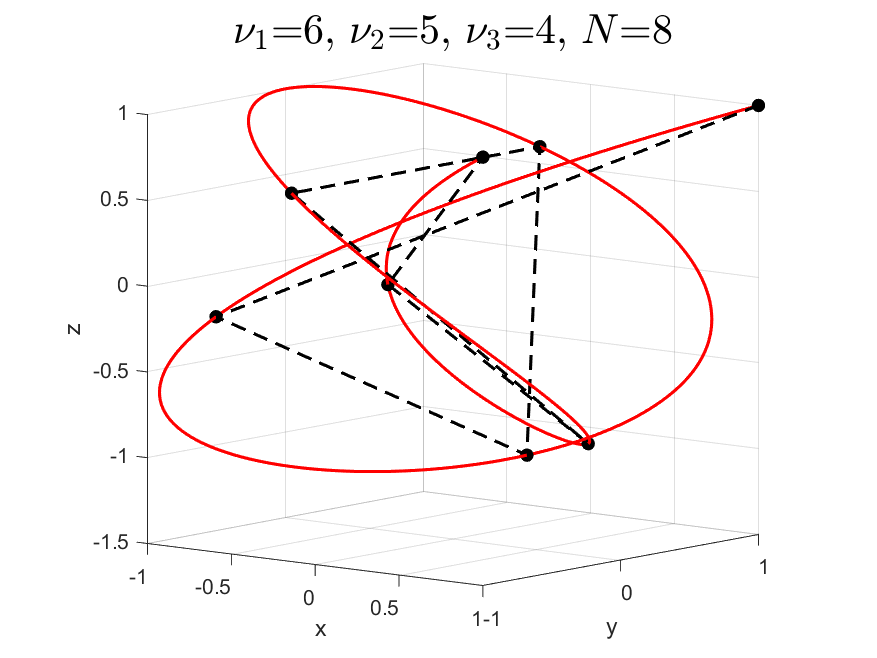}\hspace{-0.5cm}
\includegraphics[scale=0.28]{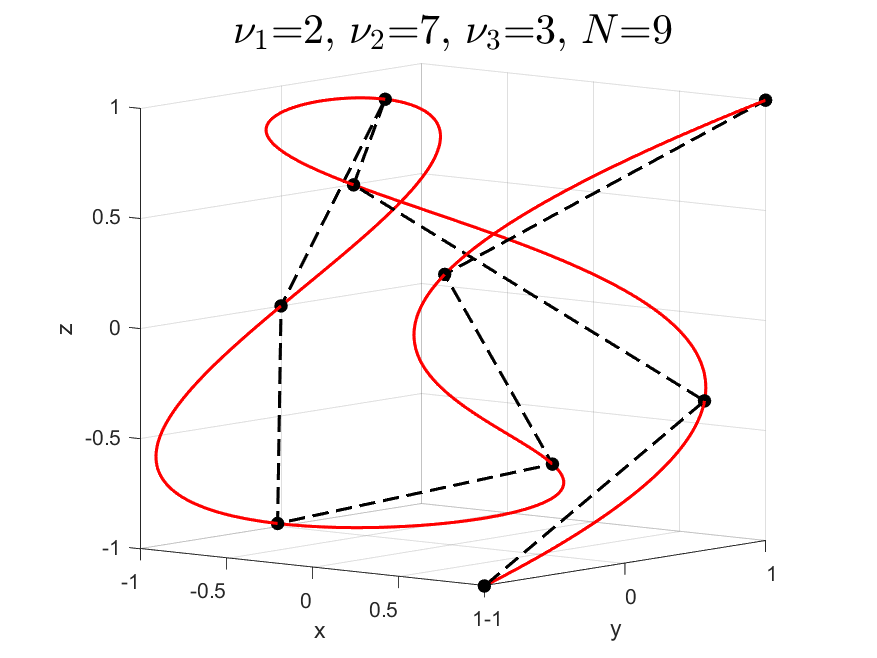}
\caption{Examples of spatial Lissajous curves with parametric equations $x(u)=\cos(\nu_1 u)$,  
$y(u)=\cos(\nu_2 u)$, $z(u)=\cos(\nu_3 u)$ with $u \in [0,\pi]$, obtained after 6 subdivision steps of the interpolatory schemes $m_{12,k}(z)$ (left), $m_{14,k}(z)$ (center), $m_{16,k}(z)$ (right) reproducing $V_{2n+2,\theta}$ with $\theta=\pi/(N-1)$ and $n=\max\{\nu_1, \nu_2, \nu_3\}$,
starting from the dashed polygons with vertices ${\bf p}_i^{(0)}=(x_i^{(0)}, y_i^{(0)}, z_i^{(0)})$, $i=1,\ldots, N$ where $x_i^{(0)}=\cos \left(\nu_1 \theta (i-1) \right)$, $y_i^{(0)}=\cos \left(\nu_2 \theta (i-1) \right)$, $z_i^{(0)}=\cos \left(\nu_3 \theta (i-1) \right)$.}
\label{fig_lissajous_cube}
\end{figure}

\begin{figure}[h!]
\centering
\includegraphics[scale=0.27]{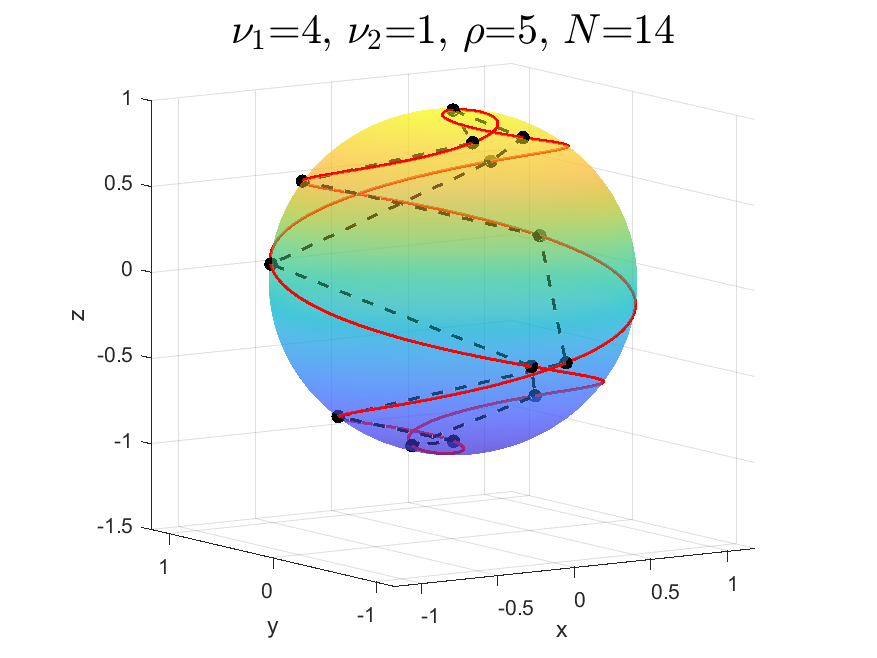}\hspace{-0.65cm}
\includegraphics[scale=0.27]{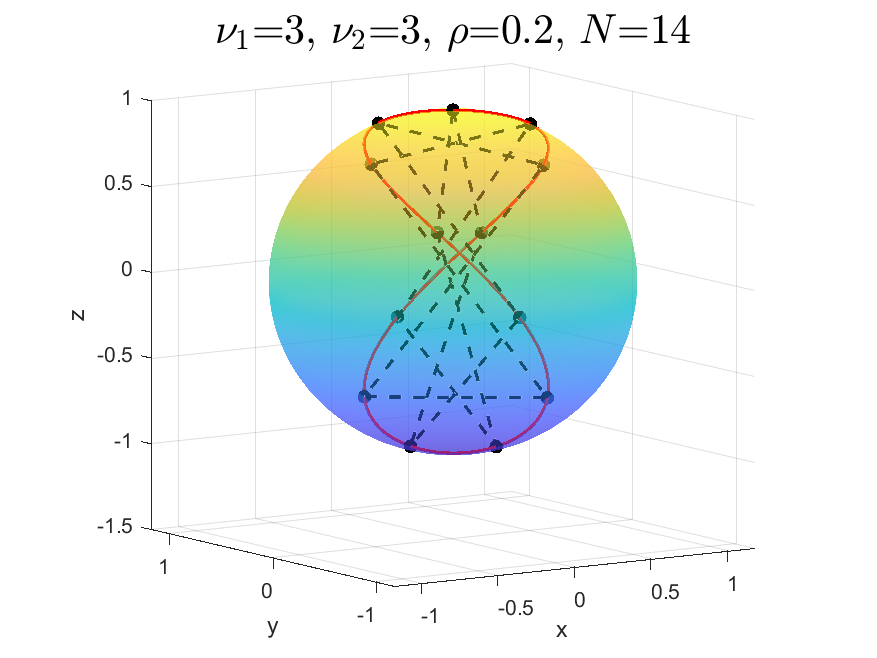}\hspace{-0.65cm}
\includegraphics[scale=0.27]{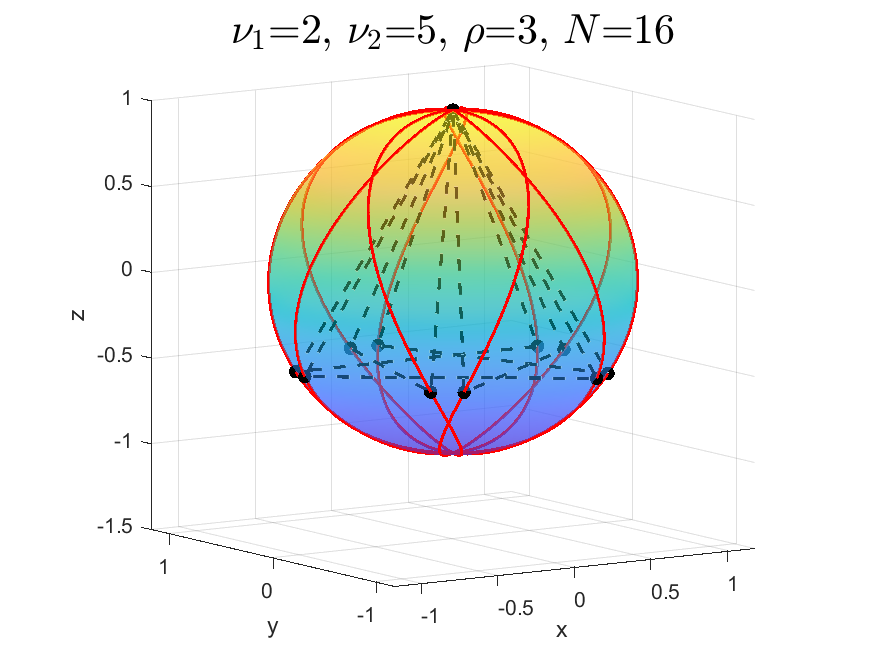}
\caption{Examples of spherical Lissajous curves with parametric equations $x(u)=\sin(\nu_2 u) \cos(\nu_1 u - \rho \pi)$,  
$y(u)=\sin(\nu_2 u) \sin(\nu_1 u - \rho \pi)$, $z(u)=\cos(\nu_2 u)$ with $u \in [0,2\pi]$, obtained after 6 subdivision steps of the interpolatory schemes $m_{12,k}(z)$ (left), $m_{14,k}(z)$ (center), $m_{16,k}(z)$ (right) reproducing $V_{2n+2,\theta}$ with $\theta=(2 \pi)/(N-1)$ and $n=\nu_1+\nu_2$, starting from the dashed polygons with vertices ${\bf p}_i^{(0)}=(x_i^{(0)}, y_i^{(0)}, z_i^{(0)})$, $i=1,\ldots, N$ where $x_i^{(0)}=\sin \left(\nu_2 \theta (i-1) \right) \cos \left(\nu_1 \theta (i-1) - \rho \pi \right)$, $y_i^{(0)}=\sin \left(\nu_2 \theta (i-1) \right) \sin \left(\nu_1 \theta (i-1) - \rho \pi \right)$, $z_i^{(0)}=\cos \left(\nu_2 \theta (i-1) \right)$.}
\label{fig_lissajous_sphere}
\end{figure}

\smallskip  
Figures  \ref{fig_starshape2d3d}, \ref{fig_lissajous}, \ref{fig_lissajous_cube}, \ref{fig_lissajous_sphere} show examples of planar and spatial curves obtained by means of  minimum-support interpolatory subdivision schemes reproducing different sets of integer powers of exponentials with values of $n$ between $4$ and $10$. 

The following proposition offers a new proof (alternative to the one in \cite{DLL03}) to show that, as $k \rightarrow +\infty$, $m_{2n+2,k}(z)$ converges to the Laurent polynomial of the $(2n+2)$-point Dubuc-Deslauriers scheme.\\

\begin{prop}\label{prop_limit}
For any $n\in\mathbb{N}$, the Laurent polynomial in \eqref{eq:a2np2} satisfies
$$
\lim_{k\rightarrow +\infty}  m_{2n+2,k}(z)\;=\; \frac{(1+z)^{2n+2}}{2^{2n+1} \, z^{n+1}} \, \sum_{s=0}^{n} \binom{n+s}{s} \, (-1)^s \, \frac{(1-z)^{2s}}{4^s z^s}.$$ 
\end{prop}

\begin{proof}
Taking into account that
$$
\displaystyle \frac{\displaystyle \prod_{\ell=0}^{i-1} (2\ell+1)}{i!} 
\; = \;
\displaystyle \frac{(2i-1)!!}{i!}  
\; = \; 
\displaystyle \frac{(2i-1)!}{2^{i-1} \, (i-1)! \, i!}  
\; = \; 2^{-i+1} \, \binom{2i-1}{i-1},
$$
the Laurent polynomials in \eqref{eq:bz_def} satisfy 
$$
\lim_{k\rightarrow +\infty} b_{i,k}(z)\;=\;(-1)^i \
2^{-2i-1} \, \binom{2i-1}{i-1} \, \ \frac{(1+z^2) \ (1-z)^{2i}}{z^{i+1}}.
$$
Thus, recalling \eqref{eq:a2np2} and  Remark \ref{rem_bspline_similarity}, we can write 
$$
\begin{array}{c}
\lim_{k\rightarrow +\infty}  m_{2n+2,k}(z) \;=\;
\frac{(1+z)^2}{2z}  \left( 1+ \frac{1+z^2}{(1+z)^2} \sum_{i=1}^n  (-1)^i \binom{2i-1}{i-1} \frac{(1-z^2)^{2i}}{2^{4i-1} z^{2i}}   \right).
\end{array}
$$
In order to prove the claimed result, we then only need to show that this limit expression coincides exactly with the one in the statement.
By introducing the notation
\[
    \mathcal{F}(n)\;:=\;\sum_{i=1}^n  (-1)^i \
 \, \binom{2i-1}{i-1} \, \frac{(1-z^2)^{2i}}{2^{4i-1} z^{2i}},
\]
and
\[
    \mathcal{G}(n)\;:=\;\frac{(1+z)^2}{1+z^2} \, \Bigg(  \frac{(1+z)^{2n}}{2^{2n} z^n} \  \sum_{s=0}^{n} \binom{n+s}{s} \, (-1)^s \, \frac{(1-z)^{2s}}{4^s z^s} -1 \Bigg),
\]
this can be equivalently done by showing that  
$$
\frac{1+z^2}{2z}\mathcal{F}(n)\;+\;\frac{(1+z)^2}{2z} \; = \; \frac{1+z^2}{2z}\mathcal{G}(n)\;+\;\frac{(1+z)^2}{2z}.
$$ 
Thus, we proceed proving by induction on $n$ that
$$ 
\begin{array}{l}
\mathcal{F}(n)\;=\;\mathcal{G}(n).
 \end{array}
$$
First we show that the statement holds for $n = 1$. Indeed,
$\mathcal{F}(1)=- (1-z^2)^{2}/(8 z^{2})$
and
$$
\begin{array}{rcl}
\mathcal{G}(1)&=&\frac{(1+z)^2}{1+z^2} \Bigg(  \frac{(1+z)^{2}}{4 z}   
\Big( 1 - \frac{(1-z)^{2}}{2 z} \Big ) -1 \Bigg)\;=\\ \\
&=&\frac{(1+z)^2}{1+z^2} \Bigg(     
  \frac{(1+z)^{2} (-1+4z-z^{2})}{8 z^2}   -1 \Bigg)
  \;=\;- \frac{(1+z)^2 (1-z)^2}{8 z^2}.
\end{array}
$$
Then, to prove the induction step we use the well-known binomial relationships
\begin{equation}\label{eq:binom1}
\binom{n-1+s}{s}\;=\;\binom{n+s}{s} \;-\; \binom{n-1+s}{s-1}    
\end{equation}
and
\begin{equation}\label{eq:binom2}
 \binom{2n-1}{n-1}\;=\;\binom{2n-1}{n},   
\end{equation}
as well as the fact that
\begin{equation}\label{eq:binom3} 
\begin{array}{l}
\binom{2n-1}{n} 
+  \binom{2n}{n} \, \frac{(1-z)^2}{4z}  
- \binom{2n-1}{n-1} \, \frac{1+z^2}{2z} 
 \;=\\ \\
\qquad\qquad\qquad
\overset{\mathrm{\eqref{eq:binom2}}}{=}
\;\binom{2n}{n} \, \frac{(1-z)^2}{4z} \;-\; \binom{2n-1}{n} \frac{(1-z)^2}{2z}
 \;=\; 0. 
 \end{array}
\end{equation}
Supposing $\mathcal{F}(j)=\mathcal{G}(j)$, for every $j\in\{1,\ldots,n-1\}$, we have
\begin{equation} \label{eq:step_1}
\begin{array}{rcl}
\mathcal{F}(n)&=& (-1)^n  \, \binom{2n-1}{n-1} \, \frac{(1-z^2)^{2n}}{2^{4n-1} z^{2n}}
\;+\;
\mathcal{F}(n-1)\\ \\
&=&
(-1)^n  \, \binom{2n-1}{n-1} \, (1+z)^{2n}\frac{(1-z)^{2n}}{2^{4n-1} z^{2n}}
\;+\;
\mathcal{G}(n-1),
\end{array}
\end{equation}
where
\[
    \mathcal{G}(n-1)\;=\; \frac{(1+z)^{2n}}{(1+z^2)  2^{2n-2} z^{n-1}} \mathcal{H}(n-1) - \frac{(1+z)^2}{1+z^2},
\]
with
\[\begin{array}{rcl}
    \mathcal{H}(n)&:=&\sum_{s=0}^{n} \binom{n+s}{s}  (-1)^s  \frac{(1-z)^{2s}}{4^s z^s}.
\end{array}\]
Now,
\[\begin{array}{l}
    \mathcal{H}(n-1)\;=\\ \\
    \qquad\overset{\eqref{eq:binom1}}{=}\; \sum_{s=0}^{n-1} \binom{n+s}{s} \, (-1)^s \, \frac{(1-z)^{2s}}{4^s z^s} \;-\;  \sum_{s=1}^{n-1} \binom{n-1+s}{s-1} \, (-1)^s \, \frac{(1-z)^{2s}}{4^s z^s}\\ \\
    \qquad=\;
    \sum_{s=0}^{n-1} \binom{n+s}{s} \, (-1)^s \, \frac{(1-z)^{2s}}{4^s z^s} \;-\; \sum_{s=0}^{n-2} \binom{n+s}{s} \, (-1)^{s+1} \, \frac{(1-z)^{2s+2}}{4^{s+1} z^{s+1}}\\ \\
    \qquad=\;
    \frac{(1+z)^2}{4z}\mathcal{H}(n)\;-\;(-1)^n\frac{(1-z)^{2n}}{4^{n}z^n}\mathcal{I}(n),
\end{array}\]
with
\[\begin{array}{rcl}
    \mathcal{I}(n)&:=&
    \binom{2n}{n}\;-\binom{2n-1}{n-1}\;+\;\binom{2n}{n}\frac{(1-z)^2}{4z}\\ \\
    &\overset{\eqref{eq:binom1}}{=}& \binom{2n-1}{n}\;+\;\binom{2n}{n}\frac{(1-z)^2}{4z}\\ \\
    &\overset{\eqref{eq:binom3}}{=}& \binom{2n-1}{n} \, \frac{1+z^2}{2z}.
\end{array}\]
In particular,
$$
\frac{\mathcal{I}(n)}{1+z^2} = \binom{2n-1}{n-1} \, \frac{1}{2z}. 
$$
Collecting everything together in \eqref{eq:step_1}, we thus obtain
\[\begin{array}{lll}
    \mathcal{F}(n) &=&
\frac{(1+z)^2 \, (1+z)^{2n}}{(1+z^2) \, 2^{2n} z^{n}} \, \mathcal{H}(n) - \frac{(1+z)^2}{1+z^2} \smallskip \\
&+& \, (-1)^n \, \frac{(1-z^2)^{2n}}{2^{4n-2} z^{2n-1}} \, \left(   
\binom{2n-1}{n-1} \, \frac{1}{2z} - \frac{\mathcal{I}(n)}{1+z^2} 
\right)  
    \\\\
    &=&\frac{(1+z)^2}{1+z^2} \, \Bigg(  \frac{(1+z)^{2n}}{2^{2n} z^n} \  \sum_{s=0}^{n} \binom{n+s}{s} \, (-1)^s \, \frac{(1-z)^{2s}}{4^s z^s} -1 \Bigg)
    \;=\; \mathcal{G}(n),
\end{array}\]
so concluding the proof.
\end{proof} 

\section{Closing remarks}\label{sec:closure}
This paper introduced new results in the domain of orthogonal polynomials.
In particular, new identities satisfied by Chebyshev polynomials of the first kind were disclosed and new connections to q-hypergeometric functions were shown.
Among them, reciprocals of Chebyshev polynomials have been proven to be special big q-Jacobi polynomials.
Immediate benefits of these new formulas have been presented 
to obtain closed-form expressions of subdivision symbols associated to minimum-support interpolating subdivision schemes reproducing finite sets of integer powers of exponentials. Since Chebyshev polynomials and their reciprocals have interesting connections to other mathematical concepts and are valuable tools in various mathematical and engineering applications, we believe that our novel identities could be of benefit also in other fields.

\bigskip
\vskip 0.1in\noindent
{\bf Acknowledgements}.  
This research has been accomplished within RITA (Research ITalian network on Approximation) and UMI-TAA. The last two authors are members of the INdAM research group GNCS (Gruppo Nazionale Calcolo Scientifico - Istituto Nazionale di Alta Matematica).


\begin{thebibliography}{99}

\bibitem{AA85}
Andrews G.E., Askey R.: Classical orthogonal polynomials. In: Brezinski C., Draux A., Magnus A.P., Maroni P., Ronveaux A. (eds), Polynômes Orthogonaux et Applications. Lecture Notes in Mathematics, vol 1171, pp. 36-62. Springer, Berlin, Heidelberg, 1985

\bibitem{BCM19}
Beccari C.V., Casciola G., Mazure M.-L.: Design or not design? A numerical characterisation for piecewise Chebyshevian splines. Numer. Algor. 81, 1-31 (2019) 

\bibitem{BCM20}
Beccari C.V., Casciola G., Mazure M.-L.: Critical length: An alternative approach. J. Comput. Appl. Math.
 370, 112603 (2020)  

\bibitem{BCR07} 
Beccari C., Casciola G., Romani L.: A non-stationary uniform tension controlled interpolating 4-
point scheme reproducing conics. Comput. Aided Geom. Des. 24(1), 1–9 (2007) 

\bibitem{CMP04}
Carnicer J.-M., Mainar E., Peña J.-M.: Critical length for design purposes and extended Chebyshev spaces. Constr. Approx. 20, 55–71 (2004)

\bibitem{CCR} 
Conti C., Cotronei M., Romani L.: Beyond B-splines: Exponential pseudo-splines and subdivision schemes reproducing exponential polynomials. Dolomites Res. Notes Approx. 10, 31–42 (2017) 

\bibitem{CGR11} 
Conti C., Gemignani L., Romani L.: From approximating to interpolatory non-stationary
subdivision schemes with the same generation properties. Adv. Comput. Math. 35, 217-241 (2011) 

\bibitem{CGR16} 
Conti C., Gemignani L., Romani L.: Exponential pseudo-splines:
Looking beyond exponential B-splines. J. Math. Anal. Appl. 439, 32–56 (2016) 

\bibitem{CR2010} 
Conti C., Romani L.: Affine combination of B-spline subdivision masks
and its non-stationary counterparts. BIT Numer. Math.  50, 269-299 (2010) 
 
\bibitem{CR2011} 
Conti C., Romani L.:  Algebraic conditions on non-stationary subdivision symbols for
exponential polynomial reproduction. J. Comput. Appl. Math. 236, 543–556 (2011) 

\bibitem{DD} 
Deslauriers G., Dubuc S.: Symmetric iterative interpolation processes. Constr. Approx. 5, 49–68 (1989)

\bibitem{DS2006} 
Dong B., Shen Z.: Linear independence of pseudo-splines. Proc. Amer. Math. Soc. 134(9), 2685-2694 (2006)

\bibitem{DHSS} 
Dyn N., Hormann K., Sabin M.A., Shen Z.: Polynomial reproduction by symmetric subdivision schemes.
J. Approx. Theory 155(1), 28-42 (2008)

\bibitem{DLL03}
Dyn N., Levin D., Luzzatto A.: Exponential reproducing subdivision schemes. Found. Comput. Math. 3, 187–206 (2003)

\bibitem{GasperRahman} 
Gasper G., Rahman M.: Basic Hypergeometric Series, 2nd ed.; Cambridge University Press: Cambridge, UK, 2004

\bibitem{GP00} 
Gori L., Pitolli F.:  A class of totally positive refinable functions. Rend. Mat. Ser. VII, 20, 305–322 (2000)

\bibitem{GIM92} 
Gupta D.P., Ismail M.E.H, Masson D.R.: 	
Contiguous relations, basic hypergeometric functions, and orthogonal polynomials. II. Associated big q-Jacobi polynomials.
J. Math. Anal. Appl. 171(2), 477-497 (1992)

\bibitem{Yoon13}
Jeong B., Kim H.O., Lee Y.J., Yoon J.: Exponential polynomial reproducing property
of non-stationary symmetric subdivision schemes and normalized exponential B-splines.
Adv. Comput. Math. 38, 647–666 (2013) 

\bibitem{Jeong13}
Jeong B., Lee Y.J., Yoon J.: A family of non-stationary subdivision schemes reproducing exponential polynomials.
J. Math. Anal. Appl. 402(1), 207-219 (2013) 

\bibitem{KSbook}
Koekoek R., Lesky P.A., Swarttouw R.F.: Hypergeometric Orthogonal Polynomials
and Their q-Analogues.
Springer-Verlag, Berlin, Heidelberg, 2010

\bibitem{MP10}  
Mainar E., Peña J.M.: Optimal bases for a class of mixed spaces and their associated spline spaces. Computers and Mathematics with Applications 59, 1509–1523
(2010) 

\bibitem{NR15} 
Novara P., Romani L.: Building blocks for designing arbitrarily smooth subdivision
schemes with conic precision. J. Comput. Appl. Math. 279, 67–79 (2015) 

\bibitem{Romani09} 
Romani L.: From approximating subdivision schemes for exponential splines to
high-performance interpolating algorithms. J. Comput. Appl. Math. 224, 383-396 (2009) 

\bibitem{RVJima} 
Romani L., Hernández Mederos V., Estrada Sarlabous J.: Exact evaluation of a class of nonstationary approximating subdivision
algorithms and related applications. IMA J. Numer. Anal. 36, 380–399 (2016) 

\end{thebibliography}
\end{document}